\documentclass[15pt]{article}
\newcommand{\red}[1]{\textcolor{red}{#1}}

 \newcommand{\pend}{\hfill \thicklines \framebox(5.5,5.5)[l]{}}
 \newenvironment{proof}{\noindent {\sc  Proof.} \rm}{\pend}
\newtheorem{remark}{Remark}[section]

\usepackage{fancyhdr,graphicx}
\usepackage{amsfonts}
\usepackage{epstopdf}
\usepackage{amsmath,bm}
\usepackage{epic,curves}
\usepackage{eepic}
\usepackage{eucal}
\usepackage{multirow} 
\usepackage{xcolor}
\usepackage{color}
\usepackage{cite}
\usepackage{amsmath}
\usepackage{multirow}
\usepackage{booktabs}
\usepackage{float}
\usepackage{algorithm}
\usepackage{algorithmic}
\usepackage{subfigure}

\newtheorem{theorem}{Theorem}[section]

\begin{document}
\title{\bf Equilibrium customer and socially optimal balking strategies in a constant retrial queue with multiple vacations and $N$-policy}
\author{Zhen Wang$^{a}$, Liwei Liu$^{a,*}$, Yiqiang Q. Zhao$^{b}$ \\
{\small\em  $^{a}$School of Science, Nanjing University of Science and Technology, Nanjing 210094, Jiangsu, China}\\
{\small\em  $^{b}$School of Mathematics and Statistics, Carleton University, 1125 Colonel By Drive, Ottawa,
ON K1S 5B6, Canada}}
\date{}
\footnotetext [1] {E-mail addresses:
NJUSTliuliwei@163.com
}
\renewcommand{\thefootnote}{*}
\footnotetext [1] {Corresponding author.}
\maketitle

\begin{abstract}
In this paper,  equilibrium strategies and optimal balking strategies of customers in a constant retrial queue with multiple vacations and the $N$-policy under two information levels, respectively, are investigated.
We assume that there is no waiting area in front of the server and an arriving customer is served immediately if the server is idle; otherwise (the server is either busy or on a vacation) it has to leave the system to join a virtual retrial orbit waiting for retrials according to the FCFS rules. After a service completion, if the system is not empty, the server becomes idle, available for serving the next customer, either a new arrival or a retried customer from the virtual retrial orbit; otherwise (if the system is empty), the server starts a vacation. Upon the completion of a vacation, the server is reactivated only if it finds at least $N$ customers in the virtual orbit; otherwise, the server continues another vacation. We study this model at two levels of information, respectively. For each level of information, we obtain both equilibrium and optimal balking strategies of customers, and make corresponding numerical comparisons. Through Particle Swarm Optimization (PSO) algorithm, we explore the impact of parameters on the equilibrium and social optimal thresholds, and obtain the trend in changes, as a function of system parameters, for the optimal social welfare, which provides guiding significance for social planners. Finally, by comparing the social welfare under two information levels, we find that whether the system information should be disclosed to customers depends on how to maintain the growth of social welfare.


\end{abstract}

\section{Introduction}
\label{sec:1}
In many service and electronic commerce systems, there exists a new trend to study the behavior of customers in queuing models. In these models, customers can decide whether to join or balk, according to a natural tendency to maximize their personal utility. To this end, from the perspective of game-theory, the decentralized behavior of customers in the queuing system has attracted extensive attentions in recent decades. Generally, queuing systems are divided into the observable case and the unobservable case depending on whether customers can obtain the information about the system upon arrival. The observable case was first studied by Naor \cite{Naor1969The}, who analyzed an $M/M/1$ queue model with a linear reward-cost structure, and obtained equilibrium and social optimal strategies. Subsequently, Naor's study was extensively extended,  see e.g. \cite{Edelson1975Congestion,johansen1980control,stidham1985optimal}. Specifically, Edelson and Hilderbrand \cite{Edelson1975Congestion} complemented the unobservable case to Naor's model. Chen and Frank \cite{chen2001state} generalized the model of Naor's, who assumed that customers and servers use the same discount rate to maximize their expected discount utility. Afterward, some authors studied equilibrium strategies in various invisible models with many different characteristics. The monograph of Hassin and Haviv \cite{Hassin2003To} summarized the main results of the subject under different levels of information.

The present paper aims to discuss equilibrium strategies and socially optimal balking strategies of customers in an $M/M/1$ constant retrial queue with multiple vacations and the $N$-policy. Customers' retrials are a common phenomenon in service systems and enterprise engineering. For example, arriving calls to a call center will be connected immediately if service staff is available, otherwise customers may have to retry for service after a random time. With the development of information technology, modern call centers may provide some levels of information to callers, e.g., the current number of customers waiting for service and/or expected waiting time, among other possibilities. Server vacation is another useful concept in modeling for situations, in which optimization of resources and/or reduction of cost are/is required. For vacation models, due to technical and cost (or other) reasons, the server might not be able to obtain the information about the current system capacity during the vacation, or it is impossible for the server to immediately return to work when the number of customers reaches a predetermined threshold, or the number of customers in the system is small at the end of the server vacation so that the server is reluctant to return to work, or return to work at the normal service rate. In addition, too frequent startups and changeovers on operations could lead to severe server wear and overhead on cost, the $N$-policy is usually used to solve this predicament, such as batch traffic transfer systems, controlled manufacturing systems and possible others. In the literature, there are relatively fewer papers studying queueing models with $N$-policy from the perspective of economic. Therefore, a model combining customer retrials, server multiple vacations, and the $N$-policy is of practical interest and is our focus of this paper.

As for studies on equilibrium balking strategies of customers, Burnetas and Economou \cite{burnetas2007equilibrium} considered queueing models with setup times under several information levels; Economou and Kanta \cite{economou2008equilibrium} discussed balking strategies for an observable queue with breakdowns and repairs; Liu, Ma and Li~\cite{liu2012equilibrium} explored an observable queue under single vacation policy; Ma, Liu and Li~\cite{ma2013equilibrium} presented equilibrium balking behavior under a multiple vacation policy; Sun, Li and Cheng-Guo~\cite{sun2016equilibrium} investigated equilibrium strategies and optimal balking strategies for an unobservable queue with double adaptive working vacations. Customers' equilibrium strategies for queue systems with retrials were also reported in the literature, for example, \cite{kulkarni1983game,zhang2012optimal} when balking is not allowed, and \cite{economou2011equilibrium,wang2013strategic,kumar2010single} when balking is allowed. Regarding models implemented with the $N$-policy, Guo and Hassin \cite{guo2011strategic,guo2012strategic} investigated models at two information levels with homogeneous and heterogeneous customers, respectively; Guo and Li \cite{guo2013strategic} addressed the same issue for systems, which are partially observable, such as the system capacity is observable, or the state of system is observable. Wang, Zhang and Huang~\cite{wang2017strategic} presented customers' strategic behavior and the social optimal problem in a constant retrial queue with the $N$-policy. Sun, Li and Tian~\cite{sun2017equilibrium} discussed equilibrium strategies and balking strategies with multiple vacations and the $N$-policy.

However, different from the previously mentioned literature on the $N$-policy, the present paper assumes that the system can be reactivated if and only if, upon the completion of a vacation, the server finds at least $N$ customers in the virtual orbit; otherwise, the server continues another vacation.

This paper studies equilibrium strategies and optimal balking strategies of customers in a queue with a constant retrial rate, multiple server vacations, and the $N$-policy under two information levels (the observable case and the unobservable case). In this system, there is no waiting area in front of the server and an arriving customer will be serviced immediately if the state of server is idle; otherwise (when the state of the server is busy or on vacation, it has to leave the system to join a virtual retrial orbit waiting for retries according to the FCFS rule.  After the completion of each service, the server will take a vacation if the system is empty, or becomes idle if there is at least one customer in the orbit.  The idle server will serve the next customer, either a new arrival or a retried customer, whichever comes earlier. The server be reactivated, upon return from the vacation, if at least $N$ customers are presented in the virtual orbit; otherwise, the server will start another vacation. For each type of information level, we determine equilibrium strategies and optimal balking strategies of customers and social welfare. For the observable case, in order to ensure that the server can be reactivated, we derive the optimal balking threshold of customers in the vacation state, which must be greater than the optimal threshold in busy state, and also greater than $N-1$. Therefore, there are three different queuing cases for the observable case, and we study the corresponding stationary distributions for the three queuing cases, and obtain the equilibrium social welfare per time unit. For the unobservable case, we derive the positive equilibrium arrival rate and optimal arrival rate, which are both unique. However, due to the complexity of equations involved, explicit expressions for the equilibrium balking thresholds of customers, socially optimal balking thresholds and optimal social welfare are not available in general. Hence, we use the Particle Swarm Optimization (PSO) algorithm to solve the complex analytic characteristics, by which the numerical optimal solution $({n^*}(1),{n^*}(2))$, and optimal social welfare ${U_s}({n^*}(0),{n^*}(1))$ and ${U_s}({\overline \lambda  ^*})$ are obtained. By comparing the numerical results for the two different information levels, respectively, we conclude that the customers' equilibrium behavior makes the system more congested than that under the socially optimal strategy, and whether the system information should be disclosed to customers depends on how to maintain the growth of the social welfare (i.e., potential demand arrivals). Obviously, in order to maximize the social welfare, which factor determines the level of information disclosure and when to disclose system information to customers are also crucial for the server or social planner. To conclude our main contributions made in this paper, we emphasize that to our best knowledge, a model, which combines features of retrials, multiple vacations, and the $N$-policy, has not been considered in the literature for the purpose of customers' equilibrium, and optimal balking strategies.

The remaining sections are organized as follows: In Section \ref{sec:2}, we describe the model in detail. We derive the corresponding stationary distributions for the three queueing cases, equilibrium thresholds and social benefit per time unit in Section \ref{sec:3}. Section \ref{sec:4} contributes to studies for the unobservable case, and we derive the equilibrium arrival rate and optimal arrival rate. Section \ref{sec:5} focuses on using numerical analysis to explore the theoretical findings in the previous sections, and compare the observable and unobservable cases of this model. Section \ref{sec:6} presents discussions and possible further studies.

\section{Model description}
\label{sec:2}
Consider  a single-server retrial queueing system with a constant retrial rate, multiple server vacations, and the exhaustive $N$-policy. We assume that customers arrive to the system according to a Poisson process with rate $\lambda$, served by a single server with exponential service rate $\mu$. There is no waiting area in front of the server. An arriving customer will be serviced immediately if the server is idle and leave the system immediately upon the completion of the service; otherwise (the server is either busy or on a vacation), it will join  a virtual retrial orbit according to the first-come, first-served discipline (FCFS). In practice, a customer in the orbit can be viewed as a customer on the waiting list. After the completion of a service, the server becomes idle and immediately searches for the customer from the top of the waiting list. The time of the search is a random variable,  exponentially distributed with rate $\theta$. In the search process, if a new customer arrives, the search will be immediately interrupted and the server will return to serving the arriving customer; otherwise (no arrivals during the search process),  the customer at the head of the waiting line will be served and will leave the system upon the completion of its service. After all customers in the system are served, or when the system becomes empty, the server will take a vacation of exponential amount of time $V$ with rate $\xi$. During the vacation time, the server will be not available to serve customers. Upon the completion of a vacation, the server will continue to another (independent) vacation with the same parameter if there are fewer than $N$ customers in the system; otherwise, the server will return from vacations (to idle state) and immediately start the same search process as that for the case, described above, when the server becomes idle from busy. This type of queue systems is referred to as the {M/M/1/MV} queue. Inter-arrival times of customers, service times and the times of retrials are assumed to be mutually independent.

\begin{figure*}
\centering
\includegraphics[width=0.7\textwidth]{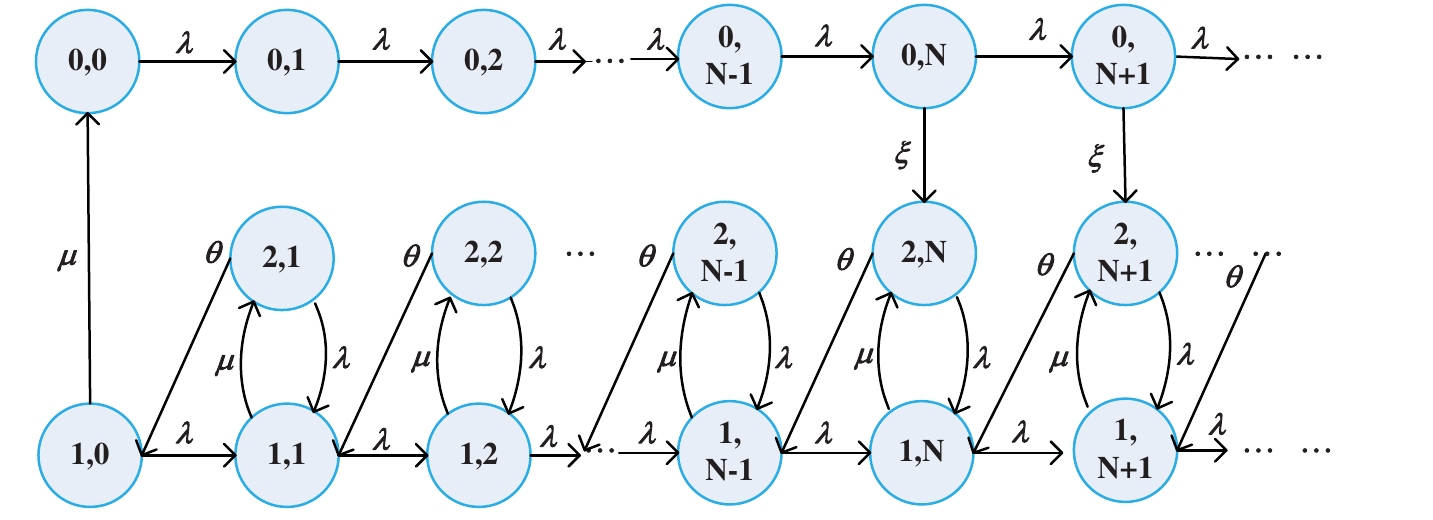}
\caption{Transition rate diagram of original model.}
\label{Fig:1}
\end{figure*}

The state of the system at time $t$ can be represented by a random vector $\{ (M(t),I(t))\} $, where $I(t)$ denotes the number of customers in the orbit, and $M(t)$ denotes the state of the server at time $t$:
\[
     M(t)=\left\{
\begin{aligned}
&0, ~~{\rm on~vacation;}\\
&1,~~{\rm busy;}\\
&2,~~{\rm idle.}
\end{aligned}
\right.
\]
Obviously, the stochastic process $\{ (M(t),I(t))\} $ is a continuous-time Markov chain. The corresponding transition rate diagram is shown in Fig.~\ref{Fig:1}. Moreover, the observable case means that arriving customers can observe all information about $M(t)$ and $I(t)$, and the unobservable case implies that arriving customers can not observe any information about $M(t)$ and $I(t)$.

The study of strategic customer behavior is important. In our case, we are interested in deciding whether an arriving customers would join or balk the system. Suppose that every customer receives the same reward of $R$ units for completing its service, which is used to quantify customer satisfaction or the added value of service. In addition, there is a waiting cost of $C$ units per time unit. The total waiting time for a customer is the continuous accumulation of the time when the customer reaches the system and until he leaves the system (including the service time). Customers are risk neutral and want to maximize their expected net benefit. Specifically, if a customer receives the reward of service is more than the expected waiting cost, then he will join the system. If a customer receives the reward of service equal to the expected waiting cost, the customer will be indifferent between joining and balking. Therefore, we only need to consider the reward with satisfying the following inequality:
\begin{equation}\label{aa1}
R > \frac{C}{\mu }.
\end{equation}
The above constrain assures that all customers, who find server being idle, always enter the system since his reward $R$ is more than the waiting cost during his expected service time ($1/\mu$). We adopt a natural linear reward-cost structure, and $U$ is defined as the expected net benefit after the service completion, i.e., $U=R-CE[W]$, where $E[W]$ is mean sojourn time of the customer. Obviously, if the customer is balking, it will generate $U=0$.

Under both levels of information levels, under the condition in (\ref{aa1}), customers will be sure to enter the system if they find that the state of the system is idle upon arrivals. However, if an arriving customer finds a vacation or busy state, he has to decide whether to leave his contact details (enter retrial orbit) or leave for ever. We further assume that the arriving customers know the policy of the system, i.e., their decisions are irrevocable: the balking customers cannot retry and customers, who joined the system,  cannot renege.

\section{The observable case}
\label{sec:3}
As mentioned above, the observable case means that arriving customers can observe all information about $M(t)$ and $I(t)$. Obviously, the information about $I(t)$ (orbit length) is useful when the arriving customer finds the server being busy or on the vacation,  but it is useless when the customer finds the server being idle upon arrivals, since in this case the customer will be served immediately regardless of the orbit length. Therefore, in all cases, the information about the system (both $M(t)$ and $I(t)$) is valuable for the arriving customer to make a better assessment on whether or not he should join the system. More specifically, if $M(t)$ is idle, the customer will join the system for sure regardless of $I(t)$; if $M(t)$ is busy or on vacation, according to the FCFS discipline, the arriving customer knows his position in the orbit, which can help him in deciding whether entering the system is preferable. In the observable case, define $W(i,n)$  to be the sojourn time of the marked customer, who joins the system at state $(i,n-1)$ ($i=1,2,3$).
For studying optimal balking strategies, we need to consider the expected (residual) net benefit of the marked customer, who is at the $n$th position in the orbit and the state of the server is $i$, after he receives the service. We denote the equilibrium balking threshold of customers and the integrated strategy at sate $i$ by ${n_e}(i)$ and $({n_e}(0),{n_e}(1))$, respectively. In addition, we denote the socially optimal balking threshold of customers and the integrated strategy at state $i$ by ${n^*}(i)$  and $({n^*}(0),{n^*}(1))$, respectively. To characterize ${n_e}(0)$ and ${n_e}(1)$, we first give the following theorem.

\begin{theorem}  \label{t31}
For the {M/M/1/MV} constant retrial queue with multiple vacations and $N$-policy, when a marked customer is at $n$th position in the orbit and the state of server is $i$ ($i=0,1,2$), the mean (residual) sojourn time $T(i,n)$ of the marked customer are given by, respectively,
\begin{equation}\label{bb1}
T(0,n) = \frac{1}{\xi } + n \cdot \frac{{\lambda  + \theta  + \mu }}{{\mu \theta }}, \quad n \ge N.
\end{equation}
\begin{equation}\label{bb2}
T(1,n) = n \cdot \frac{{\lambda  + \theta  + \mu }}{{\mu \theta }} + \frac{1}{\mu }, \quad n = 0,1 \ldots .
\end{equation}
\begin{equation}\label{bb3}
T(2,n) = n \cdot \frac{{\lambda  + \theta  + \mu }}{{\mu \theta }}, \quad n = 1,2 \ldots .
\end{equation}
\end{theorem}

\proof
Consider a marked customer arrived to the system, who found that the server is busy or on vacation. Clearly, the mean overall sojourn time of the marked customer is not affected by customers who arrive after the marked customer by finding the server being busy or on vacation, but it is affected by the customers, who enter the system after the marked customer by finding the server being idle, since in this case, by our imposed condition (\ref{aa1})) they will join the system to receive the service immediately.

Since $T(1,0)$ represents the mean residual service time of the customer, who is receiving the service, we have
\begin{equation}\label{bb4}
T(1,0) = \frac{1}{\mu }.
\end{equation}
For $n \geq 1$, let  $m(n)$ be the probability of joining the virtual orbit for the arriving customer, who finds the server being busy and $n$  customers being in the orbit. Then, based on a first step argument and noticing that the mean time to the next event is $1/(\lambda m(n) + \mu)$, and the next event is an arrival or a service completion with probability $\lambda m(n)/(\lambda m(n) + \mu)$ or $\mu/(\lambda m(n) + \mu)$, respectively, we have
\begin{equation}\label{bb5}
T(1,n) = \frac{1}{{\lambda m(n) + \mu }} + \frac{{\lambda m(n)}}{{\lambda m(n) + \mu }}T(1,n) + \frac{\mu }{{\lambda m(n) + \mu }}T(2,n), \quad n = 1, 2,  \ldots.
\end{equation}
When the server state is idle, we can similarly get
\begin{equation}\label{bb6}
T(2,n) = \frac{1}{{\lambda  + \theta }} + \frac{\lambda }{{\lambda  + \theta }}T(1,n) + \frac{\theta }{{\lambda  + \theta }}T(1,n - 1), \quad n = 1,2 \ldots.
\end{equation}
For $i=0$, we only need expressions for $n \geq N$ (see Remark~\ref{remark:1}), which is given by
\begin{equation}\label{bb7}
T(0,n) = \frac{1}{\xi } + T(2,n).
\end{equation}

In terms of (\ref{bb6}) and by solving (\ref{bb5}) for $T(1,n)$, we get
\begin{equation}\label{bb8}
T(1,n) = \frac{{\lambda  + \theta  + \mu }}{{\mu \theta }} + T(1,n - 1), \quad n = 1,2 \ldots,
\end{equation}
which leads to (\ref{bb2}). Substituting (\ref{bb2}) into (\ref{bb6}) produces (\ref{bb3}). Finally, substituting (\ref{bb3}) into (\ref{bb7}) gives (\ref{bb1}), which completes the proof. \pend

\begin{remark} \label{remark:1}
In the above theorem, we did not provide the expression for $T(0,n)$ when $n < N$, since for our purpose, we only need the expression when the server can be reactivated, or $n \geq N$.
\end{remark}

\subsection{Equilibrium}

We first study the equilibrium balking behavior of customers in the observable case, i.e., the customers can observe both information of $M(t)$ and $I(t)$ at time $t$. As mentioned above, the condition (\ref{aa1}) ensures that the customers who find the server is idle always enter the system. The customers who find a vacation or busy state have to decide whether to leave their contact details (enter retrial orbit) or leave for ever. Therefore, we only need to consider that the system is in the state of vacation or busy upon the customers arrivals.

From Theorem \ref{t31}, it is easy to know that the sojourn time $W(0,n)$ of the marked customer satisfies the following equation when he encounters the system state $(0,n)$:
\begin{equation}\label{bb12}
    E[W(0,n)] = T(0,n) = \frac{1}{\xi } + n \cdot \frac{{\lambda  + \theta  + \mu }}{{\mu \theta }}, \quad n \ge N.
\end{equation}
Define ${U_e}(0,n)=R - CE[W(0,n)]$ to be the corresponding residual net benefit of the marked customer, and solve ${U_e}(0,n) = 0$ to get the equilibrium balking threshold:
\begin{equation}\label{bb13}
{n_e}(0) = \left\lfloor {\frac{{\mu \theta }}{{\lambda  + \theta  + \mu }}(\frac{R}{C} - \frac{1}{\xi })} \right\rfloor,
\end{equation}
where the floor function $\lfloor x \rfloor$ is the largest integer smaller than $x$.

Similarly, when the server state is $i=1$, we have
\begin{equation}\label{bb14}
    E[W(1,n)] = T(1,n) =\frac{1}{\mu }+ n \cdot \frac{{\lambda  + \theta  + \mu }}{{\mu \theta }}, \quad n = 0,1 \ldots .
\end{equation}
Define ${U_e}(1,n)= R - CE[W(1,n)]$ to be the corresponding the residual net benefit of the marking customer, and slove ${U_e}(1,n) =0$ to get the balking threshold:
\begin{equation}\label{bb15}
    {n_e}(1) = \left\lfloor {\frac{{\mu \theta }}{{\lambda  + \theta  + \mu }}(\frac{R}{C} - \frac{1}{\mu })} \right\rfloor .
\end{equation}
Obviously, there are two possibilities: (i) $\mu > \xi$, which implies ${n_e}(1)>{n_e}(0)$; and (ii) $\xi > \mu$, which implies ${n_e}(0)>{n_e}(1)$. Hence, we need to discuss the stationary distribution of the system in three cases: Case~1: $N - 1 \le n(0) \le n(1)$; Case~2: $N - 1 \le n(1) \le n(0)$; and Case~3: $n(1) < N - 1 \le n(0)$ for the unobservable case. Our focus in this section is to characterize the integrated balking threshold strategy
$(n(0), n(1))$ for these three cases.

\begin{figure*}
\centering
\includegraphics[width=1.0\textwidth]{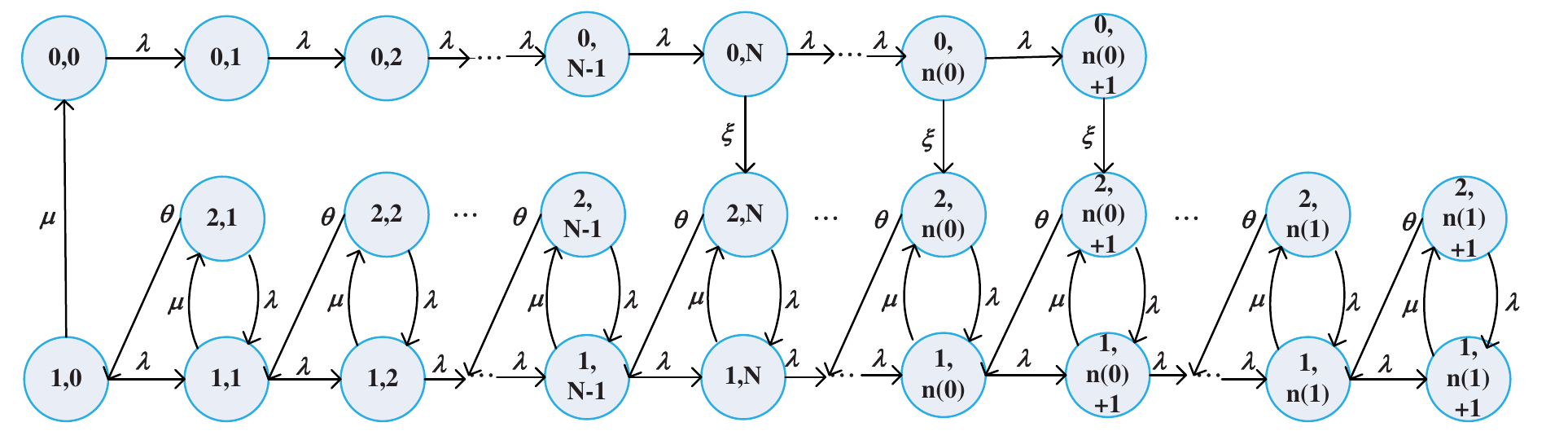}
\caption{Transition rate diagram of {($M(t)$,$I(t)$)} for the observable queues when $N-1\leq n(0) \leq n(1)$.}
\label{Fig:2}
\end{figure*}

Case~1: For $N - 1 \le n(0) \le n(1)$, the corresponding transition rate diagram is showed in Fig. \ref{Fig:2}, and the state space of $\{ (M(t),I(t))\} $ is given by:
\begin{equation} \label{eqn:case1}
    \Omega _{ob1}^e = \{ (0,n):0 \le n \le n(0) + 1\}  \cup \{ (1,n):0 \le n \le n(1) + 1\}  \cup \{ (2,n):1 \le n \le n(1) + 1\} .
\end{equation}
Define the stationary distribution as
\[
    {\pi _{i,n}} = P\{ M = i,I = n\}  = \mathop {\lim }\limits_{t \to \infty } P\{ M(t) = i,I(t) = n\} ,(i,n) \in \Omega _{ob1}^e, \quad i = 0,1,2.
\]

The stationary distribution for this case (Caese~1) is given in Theorem \ref{t32}.
\begin{theorem}  \label{t32}
For the fully observable {M/M/1/MV} constant retrial queue with multiple vacations and the $N$-policy, if $N - 1 \le n(0) \le n(1)$, then the state space $\Omega _{ob1}^e$ of $\{ (M(t),I(t))\} $ is given by (\ref{eqn:case1}), and the stationary distribution $\{ {\pi _{i,n}}\left| {(i,n) \in \Omega _{ob1}^e} \right.\} $ is given by:

\begin{equation}\label{bb16}
{\pi _{0,n}} =
\begin{cases}
\frac{\mu }{\lambda } \cdot {\pi _{1,0}}, &0 \le n \le N - 1;\\
\frac{\mu }{\lambda } \cdot {\left( {\frac{\lambda }{{\lambda  + \xi }}} \right)^{n - N + 1}} \cdot {\pi _{1,0}}, &N \le n \le n(0);\\
\frac{\mu }{\xi } \cdot {\left( {\frac{\lambda }{{\lambda  + \xi }}} \right)^{n(0) - N + 1}} \cdot {\pi _{1,0}}, &n = n(0) + 1;
\end{cases}
\end{equation}

\begin{equation}\label{bb17}
{\pi _{1,n}} =
\begin{cases}
{A_1} + {A_2} \cdot {F^n}, &0 \le n \le N - 1;\\
{B_1} + {B_2} \cdot {F^n} + {D_1} \cdot {\left( {\frac{\lambda }{{\lambda  + \xi }}} \right)^n}, &N \le n \le n(0);\\
{B_1} + {B_2} \cdot {F^{n(0) + 1}} + {D_1} \cdot \left( {\frac{{\lambda  - F\xi }}{{\lambda  + \xi }}} \right) \cdot {\left( {\frac{\lambda }{{\lambda  + \xi }}} \right)^{n(0) - 1}}\\
~~~~- (1 + \frac{\xi }{\theta }) \cdot {\left( {\frac{\lambda }{{\lambda  + \xi }}} \right)^{n(0) - N + 1}} \cdot {\pi _{1,0}}, &n = n(0) + 1;\\
{B_1} + {B_2} \cdot {F^{n(0) + 2}} + {D_1} \cdot \left( {\frac{{\lambda  - F\xi  - {F^2}\xi }}{{\lambda  + \xi }}} \right) \cdot {\left( {\frac{\lambda }{{\lambda  + \xi }}} \right)^{n(0) - 1}}\\
~~~~ - \frac{{\lambda  + \theta  + \xi  + F\theta  + F\xi }}{\theta } \cdot {\left( {\frac{\lambda }{{\lambda  + \xi }}} \right)^{n(0) - N + 1}} \cdot {\pi _{1,0}}, &n = n(0) + 2;\\
{E_2} \cdot {F^n}, &n(0) + 3 \le n \le n(1) + 1;
\end{cases}
\end{equation}
and
\begin{equation}\label{bb18}
{\pi _{2,n}} =
\begin{cases}
\frac{\mu }{{\lambda  + \theta }} \cdot {\pi _{1,n}}, &1 \le n \le N - 1;\\
\frac{\mu }{{\lambda  + \theta }} \cdot {\pi _{1,n}} + \frac{\lambda }{{\lambda  + \xi }} \cdot {\pi _{0,n}}, &N \le n \le n(0) + 1;\\
\frac{\mu }{{\lambda  + \theta }} \cdot {\pi _{1,n}}, &n(0) + 2 \le n \le n(1) + 1;
\end{cases}
\end{equation}
where
\begin{equation}\label{bb19}
\left\{ \begin{array}{l}
{A_1} = \frac{{\mu F}}{{\lambda (1 - F)}} \cdot {\pi _{1,0}},\\
{A_2} = \left( {1 + \frac{{\mu F}}{{\lambda (F - 1)}}} \right) \cdot {\pi _{1,0}},
\end{array} \right.
\end{equation}

\begin{equation}\label{bb20}
\left\{ \begin{array}{l}
{B_1} = {A_1} - {D_1} \cdot {\left( {\frac{\lambda }{{\lambda  + \xi }}} \right)^{N - 1}} \cdot \frac{{\lambda (1 - F) - F\xi }}{{(1 - F)(\lambda  + \xi )}},\\
{B_2} = {A_2} - {D_1} \cdot {\left( {\frac{\lambda }{{F(\lambda  + \xi )}}} \right)^{N - 1}} \cdot \frac{\xi }{{(1 - F)(\lambda  + \xi )}},
\end{array} \right.
\end{equation}

\begin{equation}\label{bb21}
F= \frac{{\lambda (\lambda  + \theta )}}{{\theta \mu }},
\end{equation}

\begin{equation}\label{bb22}
 {D_1} = \frac{{\xi \mu (\lambda  + \xi  + \theta ){{(\lambda  + \xi )}^{N - 1}} \cdot {\pi _{1,0}}}}{{{\lambda ^{N - 1}}\left( {(\lambda (\lambda  + \theta ) + \theta \mu ) - \lambda \theta \mu  - (\lambda  + \theta ){{(\lambda  + \xi )}^2}} \right)}},
\end{equation}
and
\begin{equation}\label{bb23}
{E_2} = \frac{{{\pi _{1,n(0) + 2}}}}{{{F^{n(0) + 2}}}},
\end{equation}
${\pi _{1,0}}$ can be obtained by the normalization condition $\sum\limits_{(i,n) \in \Omega _{ob1}^e} {{\pi _{i,n}} = 1} $.
\end{theorem}

\proof From Fig. \ref{Fig:2}, the corresponding balance equations of the stationary distribution are given as follows,
\begin{align} \label{bb24}
    \lambda {\pi _{0,0}}  &=  \mu {\pi _{1,0}}, \\
\label{bb25}
    \lambda {\pi _{0,n}} &= \lambda {\pi _{0,n - 1}}, \qquad 1 \le n \le N - 1, \\
\label{bb26}
    (\lambda  + \xi ){\pi _{0,n}} &= \lambda {\pi _{0,n - 1}}, \qquad N \le n \le n(0), \\
\label{bb27}
    \xi {\pi _{0,n(0) + 1}} &= \lambda {\pi _{0,n(0)}}, \\
\label{bb28}
    (\lambda  + \mu ){\pi _{1,0}} &= \theta {\pi _{2,1}}, \\
\label{bb29}
    (\lambda  + \mu ){\pi _{1,n}} &= \lambda {\pi _{1,n - 1}} + \lambda {\pi _{2,n}} + \theta {\pi _{2,n + 1}}, \qquad 1 \le n \le n(1), \\
\label{bb30}
    \mu {\pi _{1,n(1) + 1}} &= \lambda {\pi _{1,n(1)}} + \lambda {\pi _{2,n(1) + 1}}, \\
\label{bb31}
    (\lambda  + \theta ){\pi _{2,n}} &= \mu {\pi _{1,n}}, \qquad 1 \le n \le N - 1 \text{ and } n(0) + 2 \le n \le n(1) + 1, \\
\label{bb32}
    (\lambda  + \theta ){\pi _{2,n}} &= \mu {\pi _{1,n}} + \xi {\pi _{0,n}}, \qquad N \le n \le n(0) + 1.
\end{align}
We first consider the stationary distribution $\{ {\pi _{0,n}}\left| {0 \le n \le n(0) + 1} \right.\} $. From (\ref{bb24}) and (\ref{bb25}), we can obtain
\begin{equation}\label{bb33}
{\pi _{0,n}} = \frac{\mu }{\lambda }{\pi _{1,0}}, \quad 0 \le n \le N - 1.
\end{equation}
From (\ref{bb26}) and (\ref{bb33}),
\begin{equation}\label{bb34}
\begin{array}{l}
{\pi _{0,n}} = \frac{\lambda }{{\lambda  + \xi }}{\pi _{0,n - 1}} = {\left( {\frac{\lambda }{{\lambda  + \xi }}} \right)^{n - N + 1}}{\pi _{0,N - 1}}\\
{\rm{                          }} ~~~~~~~~~~~~~~~~~~~~~~~= {\left( {\frac{\lambda }{{\lambda  + \xi }}} \right)^{n - N + 1}}\frac{\mu }{\lambda }{\pi _{1,0}}\\
{\rm{                          }} ~~~~~~~~~~~~~~~~~~~~~~~= \frac{\mu }{\lambda }{\left( {\frac{\lambda }{{\lambda  + \xi }}} \right)^{n - N + 1}}{\pi _{1,0}}, \qquad N \le n \le n(0).
\end{array}
\end{equation}
Based on (\ref{bb27}) and (\ref{bb34}), we can get
\begin{equation}\label{bb35}
{\pi _{0,n(0) + 1}} = \frac{\lambda }{\xi }{\left( {\frac{\lambda }{{\lambda  + \xi }}} \right)^{n(0) - N + 1}}\frac{\mu }{\lambda }{\pi _{1,0}} = \frac{\mu }{\xi }{\left( {\frac{\lambda }{{\lambda  + \xi }}} \right)^{n(0) - N + 1}}{\pi _{1,0}}.
\end{equation}
Therefore, we can get (\ref{bb16}) from the above discussion.

Next, we consider the stationary distribution $\{ {\pi _{1,n}}\left| {0 \le n \le N-1} \right.\} $. From (\ref{bb29}) and (\ref{bb31}), we can obtain
\begin{equation}\label{bb36}
(\lambda  + \mu ){\pi _{1,n}} = \lambda {\pi _{1,n - 1}} + \frac{{\lambda \mu }}{{\lambda  + \theta }}{\pi _{1,n}} + \frac{{\theta \mu }}{{\lambda  + \theta }}{\pi _{1,n + 1}},1 \le n \le N - 1.
\end{equation}
The solution of (\ref{bb36}) can be given by the following homogeneous linear difference equation:
\begin{equation}\label{bb37}
\frac{{\theta \mu }}{{\lambda  + \theta }}{x_{n + 1}} - \left( {\lambda  + \frac{{\theta \mu }}{{\lambda  + \theta }}} \right){x_n} + \lambda {x_{n - 1}} = 0,1 \le n \le N - 1.
\end{equation}
The characteristic equation corresponding to (\ref{bb37}) is
\begin{equation}\label{bb38}
\frac{{\theta \mu }}{{\lambda  + \theta }}{x^2} - \left( {\lambda  + \frac{{\theta \mu }}{{\lambda  + \theta }}} \right)x + \lambda  = 0,
\end{equation}
which has two roots: 1 and $F = \frac{{\lambda (\lambda  + \theta )}}{{\theta \mu }}$.  Let $x_n^h = {A_1} + {A_2}{F^n}$ be the general solution of (\ref{bb37}), where ${A_1}$ and ${A_2}$ are the coefficients that need to be determined. From (\ref{bb28}) and (\ref{bb31}), we can obtain
\begin{equation}\label{bb39}
\left\{ \begin{array}{l}
{A_1} + {A_2} = {\pi _{1,0}},\\
(\lambda  + \mu )({A_1} + {A_2}) = \theta {\pi _{2,1}} = \frac{{\theta \mu }}{{\lambda  + \theta }}{\pi _{1,1}} = \frac{{\theta \mu }}{{\lambda  + \theta }}({A_1} + {A_2}F),
\end{array} \right.
\end{equation}
which yields
\begin{equation}\label{bb40}
\left\{ \begin{array}{l}
{A_1} = \frac{{\mu F}}{{\lambda (1 - F)}}{\pi _{1,0}},\\
{A_2} = \left( {1 + \frac{{\mu F}}{{\lambda (F - 1)}}} \right){\pi _{1,0}}.
\end{array} \right.
\end{equation}
Therefore,
\begin{equation}\label{bb41}
{\pi _{1,n}}={A_1} + {A_2} \cdot {F^n}, \qquad 0 \le n \le N - 1,
\end{equation}
where ${A_1}$ and ${A_2}$ are given by (\ref{bb40}).

Now, let us continue to consider the stationary distribution $\{ {\pi _{1,n}}\left| {N \le n\le \\ n(0)} \right.\} $. From (\ref{bb29}) and (\ref{bb32}), we can obtain
\begin{align}\label{bb42}
   \nonumber &  \frac{{\theta \mu }}{{\lambda  + \theta }}{\pi _{1,n + 1}} - \left( {\lambda  + \frac{{\theta \mu }}{{\lambda  + \theta }}} \right){\pi _{1,n}} + \lambda {\pi _{1,n - 1}} \\
  &  =  - \left( {1 + \frac{\theta }{{\lambda  + \xi }}} \right)\frac{{\xi \mu }}{{\lambda  + \theta }}{\left( {\frac{\lambda }{{\lambda  + \xi }}} \right)^{n - N + 1}}{\pi _{1,0}}, \qquad
 N \le n \le n(0).
\end{align}
The solutions of (\ref{bb42}) can be obtained through solving the following system of nonhomogeneous linear difference equations:
\begin{equation}\label{bb43}
\frac{{\theta \mu }}{{\lambda  + \theta }}{x_{n + 1}} - \left( {\lambda  + \frac{{\theta \mu }}{{\lambda  + \theta }}} \right){x_n} + \lambda {x_{n - 1}} =  - \left( {1 + \frac{\theta }{{\lambda  + \xi }}} \right)\frac{{\xi \mu }}{{\lambda  + \theta }}{\left( {\frac{\lambda }{{\lambda  + \xi }}} \right)^{n - N + 1}}{\pi _{1,0}}, \qquad N \le n \le n(0),
\end{equation}
whose corresponding characteristic equation is given by
\begin{equation}\label{bb44}
\frac{{\theta \mu }}{{\lambda  + \theta }}{x^2} - \left( {\lambda  + \frac{{\theta \mu }}{{\lambda  + \theta }}} \right)x + \lambda  =  - \left( {1 + \frac{\theta }{{\lambda  + \xi }}} \right)\frac{{\xi \mu }}{{\lambda  + \theta }}{\left( {\frac{\lambda }{{\lambda  + \xi }}} \right)^{n - N + 1}}{\pi _{1,0}}.
\end{equation}
Define $y_n^g = y_n^h + y_n^s$ as the general solution of (\ref{bb44}), where $y_n^h$ is the general solution of the homogeneous version of (\ref{bb44}), which is $y_n^h = {B_1} + {B_2}{F^n}$, and $y_n^s$ is a specific solution of (\ref{bb44}).

We consider a specific solution $y_n^s = {D_1}{\left( {\frac{\lambda }{{\lambda  + \xi }}} \right)^n}$ of (\ref{bb44}). Substituting it into (\ref{bb44}), we can obtain
\begin{equation}\label{bb45}
{D_1} = \frac{{\xi \mu (\lambda  + \xi  + \theta ){{(\lambda  + \xi )}^{N - 1}} \cdot {\pi _{1,0}}}}{{{\lambda ^{N - 1}}\left( {(\lambda (\lambda  + \theta ) + \theta \mu ) - \lambda \theta \mu  - (\lambda  + \theta ){{(\lambda  + \xi )}^2}} \right)}}.
\end{equation}
Thus,
\begin{equation}\label{bb46}
y_n^g = {B_1} + {B_2}{F^n} + {D_1}{\left( {\frac{\lambda }{{\lambda  + \xi }}} \right)^n}, \qquad
N \le n \le n(0),
\end{equation}
where ${B_1}$ and ${B_2}$ are the coefficients that need to be determined. By considering (\ref{bb41}), we get
\begin{equation}\label{bb47}
\left\{ \begin{array}{l}
{B_1} = {A_1} - {D_1} \cdot {\left( {\frac{\lambda }{{\lambda  + \xi }}} \right)^{N - 1}} \cdot \frac{{\lambda (1 - F) - F\xi }}{{(1 - F)(\lambda  + \xi )}},\\
{B_2} = {A_2} - {D_1} \cdot {\left( {\frac{\lambda }{{F(\lambda  + \xi )}}} \right)^{N - 1}} \cdot \frac{\xi }{{(1 - F)(\lambda  + \xi )}}.
\end{array} \right.
\end{equation}
Therefore,
\begin{equation}\label{bb48}
{\pi _{1,n}}={B_1} + {B_2} \cdot {F^n} + {D_1} \cdot {\left( {\frac{\lambda }{{\lambda  + \xi }}} \right)^n}, \qquad N \le n \le n(0),
\end{equation}
where ${D_1}$, ${B_1}$ and ${B_2}$ are given by (\ref{bb45}) and (\ref{bb47}), respectively. Specially, based on (\ref{bb29}), (\ref{bb32}), (\ref{bb16}) and (\ref{bb48}), we can obtain the stationary distribution of $\{ {\pi _{1,n(0) + 1}}\} $ as follows:
\begin{equation}\label{bb49}
{\pi _{1,n(0) + 1}}= {B_1} + {B_2} \cdot {F^{n(0) + 1}} + {D_1} \cdot \left( {\frac{{\lambda  - F\xi }}{{\lambda  + \xi }}} \right) \cdot {\left( {\frac{\lambda }{{\lambda  + \xi }}} \right)^{n(0) - 1}} - (1 + \frac{\xi }{\theta }) \cdot {\left( {\frac{\lambda }{{\lambda  + \xi }}} \right)^{n(0) - N + 1}} \cdot {\pi _{1,0}}.
\end{equation}
Based on (\ref{bb29}), (\ref{bb31}), (\ref{bb32}), (\ref{bb16}) and (\ref{bb49}), we can get the stationary distribution of $\{ {\pi _{1,n(0) + 2}}\} $ as follows:
\begin{equation}\label{bb50}
\begin{split}
{\pi _{1,n(0) + 2}} = {B_1} + {B_2} \cdot {F^{n(0) + 2}} + {D_1} \cdot \left( {\frac{{\lambda  - F\xi  - {F^2}\xi }}{{\lambda  + \xi }}} \right) \cdot {\left( {\frac{\lambda }{{\lambda  + \xi }}} \right)^{n(0) - 1}} \\
- \frac{{\lambda  + \theta  + \xi  + F\theta  + F\xi }}{\theta } \cdot {\left( {\frac{\lambda }{{\lambda  + \xi }}} \right)^{n(0) - N + 1}} \cdot {\pi _{1,0}}.
\end{split}
\end{equation}

Continue our proof for the case of $\{ {\pi _{1,n}}\left| {n(0) + 3 \le n \le n(1) + 1} \right.\} $. In this case, the general solution of (\ref{bb37}) is $z_n^h = {E_1} + {E_2}{F^n}$, where ${E_1}$ and ${E_2}$ are the coefficients that need to be determined. From (\ref{bb30}), (\ref{bb31}) and (\ref{bb50}), we can obtain ${E_1}=0$ and
\begin{equation}\label{bb51}
{E_2} = \frac{{{\pi _{1,n(0) + 2}}}}{{{F^{n(0) + 2}}}}.
\end{equation}
Therefore,
\begin{equation}\label{bb52}
{\pi _{1,n}}={E_2} \cdot {F^n},n(0) + 3 \le n \le n(1) + 1,
\end{equation}
which leads to (\ref{bb17}).

Finally, we consider the stationary distribution of $\{ {\pi _{2,n}}\left| {1 \le n \le n(1) + 1} \right.\} $. From (\ref{bb17}), (\ref{bb31}) and (\ref{bb32}), we can easily get (\ref{bb18}).

In summary, (\ref{bb16}), (\ref{bb17}) and (\ref{bb18}) are all related to ${\pi _{0,1}}$, and we can get ${\pi _{0,1}}$ by normalizing conditions $\sum\limits_{(i,n) \in \Omega _{ob1}^e} {{\pi _{i,n}} = 1} $.  \pend
\vspace*{5mm}

Based on Fig. \ref{Fig:2} and Theorem \ref{t32}, we know that the balking states of customers are $(0,n(0) + 1)$ and $(1,n(1) + 1)$. For the social optimization, which will be considered later, we define $U_{ob1}^e(n(0),n(1))$ to be the social benefit per time unit in Case~1: $N - 1 \le n(0) \le n(1)$, or
\begin{equation}\label{bb53}
U_{ob1}^e (n(0),n(1)) = \lambda R(1 - {\pi _{0,n(0) + 1}} - {\pi _{1,n(1) + 1}}) - C(\sum\limits_{n = 0}^{n(0) + 1} {n{\pi _{0,n}} + } \sum\limits_{n = 0}^{n(1) + 1} {n{\pi _{1,n}}} + \sum\limits_{n = 1}^{n(1) + 1} {n{\pi _{2,n}}}).
\end{equation}
Indeed, the first summand of (\ref{bb53}) is the effective arrival rate at the system times the reward $R$, while the second summand is the mean number of customers in the system. Obviously, the equilibrium social benefit is $U_{ob1}^e({n_e}(0),{n_e}(1))$.
\vspace*{5mm}

\begin{figure*}
\centering
\includegraphics[width=1.0\textwidth]{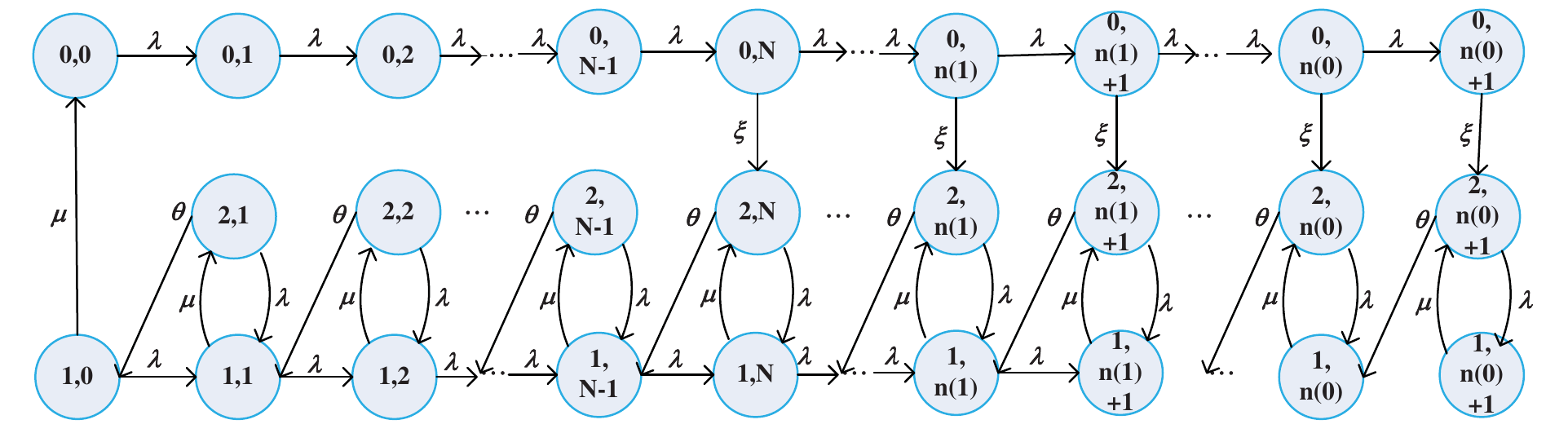}
\caption{Transition rate diagram of {($M(t)$,$I(t)$)} for the observable queues when $N-1\leq n(1) \leq n(0)$.}
\label{Fig:3}
\end{figure*}

Case~2: For $N - 1 \le n(1) \le n(0)$, the corresponding transition rate diagram is showed in Fig. \ref{Fig:3}, and the state space of $\{ (M(t),I(t))\} $ is given by
\begin{equation}\label{eqn:case2}
    \Omega _{ob2}^e = \{ (0,n):0 \le n \le n(0) + 1\}  \cup \{ (1,n):0 \le n \le n(0) + 1\}  \cup \{ (2,n):1 \le n \le n(0) + 1\} .
\end{equation}

The stationary distribution for this case (Caese~2) is given in Theorem~\ref{t33}.
\begin{theorem}  \label{t33}
For the fully observable {M/M/1/MV} constant retrial queue with multiple vacations and the $N$-policy, if $N - 1 \le n(1) \le n(0)$, then the state space $\Omega _{ob2}^e$ of $\{ (M(t),I(t))\} $ is given by (\ref{eqn:case2}), and the stationary distribution $\{ {\pi _{i,n}}\left| {(i,n) \in \Omega_{ob2}^e} \right.\} $ is given by:
\begin{equation}\label{bb54}
{\pi _{0,n}} =
\begin{cases}
\frac{\mu }{\lambda } \cdot {\pi _{1,0}}, &0 \le n \le N - 1;\\
\frac{\mu }{\lambda } \cdot {\left( {\frac{\lambda }{{\lambda  + \xi }}} \right)^{n - N + 1}} \cdot {\pi _{1,0}}, &N \le n \le n(0);\\
\frac{\mu }{\xi } \cdot {\left( {\frac{\lambda }{{\lambda  + \xi }}} \right)^{n(0) - N + 1}} \cdot {\pi _{1,0}}, &n = n(0) + 1;
\end{cases}
\end{equation}
\begin{equation}\label{bb55}
{\pi _{1,n}} =
\begin{cases}
{A_1} + {A_2} \cdot {F^n}, &0 \le n \le N - 1;\\
{B_1} + {B_2} \cdot {F^n} + {D_1} \cdot {\left( {\frac{\lambda }{{\lambda  + \xi }}} \right)^n}, &N \le n \le n(1);\\
\begin{split}
{B_1} + {B_2} \cdot {F^{n(1) + 1}} + {D_1} \cdot \left( {\frac{{\lambda  - F\xi }}{{\lambda  + \xi }}} \right) \cdot {\left( {\frac{\lambda }{{\lambda  + \xi }}} \right)^{n(1) - 1}}\\
 -(\frac{\xi }{{\lambda  + \xi }} + \frac{\xi }{\theta }) \cdot {\left( {\frac{\lambda }{{\lambda  + \xi }}} \right)^{n(1) - N + 1}} \cdot {\pi _{1,0}},
\end{split}&n = n(1) + 1;\\
(1 - F){B_1} + {D_1}\left( {\frac{\lambda }{{\lambda  + \xi }} - F} \right){\left( {\frac{\lambda }{{\lambda  + \xi }}} \right)^{n(1) - 1}} + \psi (n) \cdot {\pi _{1,0}}, &n(1) + 2 \le n \le n(0)\\
{\pi _{1,n(0)}} - \left( {1 + \frac{\xi }{\theta }} \right){\left( {\frac{\lambda }{{\lambda  + \xi }}} \right)^{n(0) - N + 1}}{\pi _{1,0}},, &n=n(0) + 1;
\end{cases}
\end{equation}
and
\begin{equation}\label{bb56}
{\pi _{2,n}} =
\begin{cases}
\frac{\mu }{{\lambda  + \theta }} \cdot {\pi _{1,n}}, &1 \le n \le N - 1;\\
\frac{\mu }{{\lambda  + \theta }} \cdot {\pi _{1,n}} + \frac{\xi }{{\lambda  + \xi }} \cdot {\pi _{0,n}}, &N \le n \le n(0) + 1;
\end{cases}
\end{equation}
where ${A_i}$ ($i = 1,2$), ${B_i}$ ($i = 1,2$), $F$ and ${D_1}$ are given by (\ref{bb19}), (\ref{bb20}), (\ref{bb21}) and (\ref{bb22}), respectively.
\begin{equation}\label{bb57}
\psi (n) = \left( {\frac{\lambda }{{\lambda  + \xi }}} \right.\left. { + \frac{\lambda }{\theta }} \right){\left( {\frac{\lambda }{{\lambda  + \xi }}} \right)^{n - N}} - \left( {1 + \frac{{\lambda  + \xi }}{\theta } + \frac{\xi }{\lambda } + \frac{{\xi (\lambda  + \xi )}}{{\lambda  + \theta }}} \right){\left( {\frac{\lambda }{{\lambda  + \xi }}} \right)^{n(1) + 2 - N}}.
\end{equation}
${\pi _{1,0}}$ can be obtained by the normalization condition $\sum\limits_{(i,n) \in \Omega _{ob2}^e} {{\pi _{i,n}} = 1} $.
\end{theorem}

{\bf Proof}~~From Fig. \ref{Fig:3}, the corresponding balance equations of the  stationary distribution are given as follows:
\begin{align}\label{bb58}
    \lambda {\pi _{0,0}} &= \mu {\pi _{1,0}}, \\
\label{bb59}
    \lambda {\pi _{0,n}} &= \lambda {\pi _{0,n - 1}}, \qquad 1 \le n \le N - 1, \\
\label{bb60}
    (\lambda  + \xi ){\pi _{0,n}} &= \lambda {\pi _{0,n - 1}}, \qquad N \le n \le n(0), \\
\label{bb61}
    \xi {\pi _{0,n(0) + 1}} &= \lambda {\pi _{0,n(0)}}, \\
\label{bb62}
    (\lambda  + \mu ){\pi _{1,0}} &= \theta {\pi _{2,1}}, \\
\label{bb63}
    (\lambda  + \mu ){\pi _{1,n}} &= \lambda {\pi _{1,n - 1}} + \lambda {\pi _{2,n}} + \theta {\pi _{2,n + 1}}, \qquad 1 \le n \le n(1), \\
\label{bb64}
    \mu {\pi _{1,n(1) + 1}} &= \lambda {\pi _{1,n(1)}} + \lambda {\pi _{2,n(1) + 1}} + \theta {\pi _{2,n(1) + 2}}, \\
\label{bb65}
    \mu {\pi _{1,n}} &= \lambda {\pi _{2,n}} + \theta {\pi _{2,n + 1}}, \qquad n(1) + 2 \le n \le n(0), \\
\label{bb66}
    \mu {\pi _{1,n(0) + 1}} &= \lambda {\pi _{2,n(0) + 1}}, \\
\label{bb67}
    (\lambda  + \theta ){\pi _{2,n}} &= \mu {\pi _{1,n}}, \qquad 1 \le n \le N - 1, \\
\label{bb68}
    (\lambda  + \theta ){\pi _{2,n}} &= \mu {\pi _{1,n}} + \xi {\pi _{0,n}}, \qquad N \le n \le n(0) + 1.
\end{align}

We first consider the stationary distribution $\{ {\pi _{0,n}}\left| {0 \le n \le n(0) + 1} \right.\}$. From  (\ref{bb58})--(\ref{bb61}), the discussion is similar to the discussion for (\ref{bb24})--(\ref{bb27}), which leads to (\ref{bb54}).

We next consider the stationary distribution $\{ {\pi _{1,n}}\left| {0 \le n \le N-1} \right.\} $. From (\ref{bb62}), (\ref{bb63}) and (\ref{bb67}), the discussion is similar to that for (\ref{bb41}), and we can obtain
\begin{equation}\label{bb69}
{\pi _{1,n}}={A_1} + {A_2} \cdot {F^n},0 \le n \le N - 1,
\end{equation}
where ${A_1}$ and ${A_2}$ are given by (\ref{bb40}).

We now continue to consider the stationary distribution $\{ {\pi _{1,n}}\left| {N \le n\le \\ n(1)} \right.\} $. From (\ref{bb54}), (\ref{bb63}) and (\ref{bb68}), the discussion is similar to that for (\ref{bb48}), and we can obtain
\begin{equation}\label{bb70}
{\pi _{1,n}}={B_1} + {B_2} \cdot {F^n} + {D_1} \cdot {\left( {\frac{\lambda }{{\lambda  + \xi }}} \right)^n}, \qquad N \le n \le n(1),
\end{equation}
where ${D_1}$, ${B_1}$ and ${B_2}$ are given by (\ref{bb45}) and (\ref{bb47}), respectively. Specially, based on (\ref{bb54}), (\ref{bb63}) and (\ref{bb70}), we can get the stationary distribution of $\{ {\pi _{1,n(1) + 1}}\} $ as follows:
\begin{equation}\label{bb71}
{\pi _{1,n(1) + 1}} = {B_1} + {B_2} \cdot {F^{n(1) + 1}} + {D_1} \cdot \left( {\frac{{\lambda  - F\xi }}{{\lambda  + \xi }}} \right) \cdot {\left( {\frac{\lambda }{{\lambda  + \xi }}} \right)^{n(1) - 1}} - (\frac{\xi }{{\lambda  + \xi }} + \frac{\xi }{\theta }) \cdot {\left( {\frac{\lambda }{{\lambda  + \xi }}} \right)^{n(1) - N + 1}} \cdot {\pi _{1,0}}.
\end{equation}
Based on (\ref{bb54}), (\ref{bb64}), (\ref{bb68}) and (\ref{bb71}), we can get the stationary distribution of $\{ {\pi _{1,n(1) + 2}}\} $ as follows
\begin{equation}\label{bb72}
\begin{split}
{\pi _{1,n(1) + 2}} = (1 - F){B_1} + {D_1}\left( {\frac{\lambda }{{\lambda  + \xi }} - F} \right){\left( {\frac{\lambda }{{\lambda  + \xi }}} \right)^{n(1) - 1}}\\
- \left( {\frac{\xi }{\lambda } + \frac{{\xi (\lambda  + \xi )}}{{\lambda  + \theta }} + \frac{\xi }{{\lambda  + \xi }} + \frac{\xi }{\theta }} \right){\left( {\frac{\lambda }{{\lambda  + \xi }}} \right)^{n(1) - N + 2}}{\pi _{1,0}}.
\end{split}
\end{equation}

Continue further to consider the stationary distribution of $\{ {\pi _{1,n}}\left| {n(1) + 3 \le n \le n(0) } \right.\} $. Based on (\ref{bb54}), (\ref{bb65}) and (\ref{bb68}), we can obtain that
\begin{equation}\label{bb73}
{\pi _{1,n}} - {\pi _{1,n - 1}} =  - \left( {\frac{\xi }{{\lambda  + \xi }} + \frac{\xi }{\theta }} \right){\left( {\frac{\lambda }{{\lambda  + \xi }}} \right)^{n - N}}{\pi _{1,0}}.
\end{equation}
Recursively using (\ref{bb73}), we can obtain that
\begin{equation}\label{bb74}
{\pi _{1,n}} = \left( {\frac{\lambda }{{\lambda  + \xi }} + \frac{\lambda }{\theta }} \right)\left( {{{\left( {\frac{\lambda }{{\lambda  + \xi }}} \right)}^{n - N}} - {{\left( {\frac{\lambda }{{\lambda  + \xi }}} \right)}^{n(1) + 2 - N}}} \right) + {\pi _{1,n(1) + 2}}, \quad n(1) + 2 \le n \le n(0).
\end{equation}
Thus,
\begin{equation}\label{bb75}
{\pi _{1,n}} = (1 - F){B_1} + {D_1}\left( {\frac{\lambda }{{\lambda  + \xi }} - F} \right){\left( {\frac{\lambda }{{\lambda  + \xi }}} \right)^{n(1) - 1}} + \psi (n){\pi _{1,0}},n(1) + 2 \le n \le n(0),
\end{equation}
where
\begin{equation}\label{bb76}
\psi (n) = \left( {\frac{\lambda }{{\lambda  + \xi }}} \right.\left. { + \frac{\lambda }{\theta }} \right){\left( {\frac{\lambda }{{\lambda  + \xi }}} \right)^{n - N}} - \left( {1 + \frac{{\lambda  + \xi }}{\theta } + \frac{\xi }{\lambda } + \frac{{\xi (\lambda  + \xi )}}{{\lambda  + \theta }}} \right){\left( {\frac{\lambda }{{\lambda  + \xi }}} \right)^{n(1) + 2 - N}}.
\end{equation}
Specially, based on (\ref{bb54}), (\ref{bb65}) and (\ref{bb68}), we can get the stationary distribution of $\{ {\pi _{1,n(0) + 1}}\} $ as follows:
\begin{equation}\label{bb77}
{\pi _{1,n(0) + 1}} = {\pi _{1,n(0)}} - \left( {1 + \frac{\xi }{\theta }} \right){\left( {\frac{\lambda }{{\lambda  + \xi }}} \right)^{n(0) - N + 1}}{\pi _{1,0}}.
\end{equation}
Therefore, taking into account all the above discussions, we can get (\ref{bb55}).

Finally, we consider the stationary distribution of $\{ {\pi _{2,n}}\left| {1 \le n \le n(0) + 1} \right.\} $. From (\ref{bb55}), (\ref{bb67}) and (\ref{bb68}), We can easily get (\ref{bb56}).

In summary, (\ref{bb54}), (\ref{bb55}) and (\ref{bb56}) are all related to ${\pi _{0,1}}$, and we can get ${\pi _{0,1}}$ by normalizing conditions $\sum\limits_{(i,n) \in \Omega _{ob2}^e} {{\pi _{i,n}} = 1} $. \pend
\vspace*{5mm}

Based on Fig \ref{Fig:3} and Theorem \ref{t33}, we know that the states at which customers will balk are $(0,n(0) + 1)$ and $(1,n(1) + 1)$. Denote $U_{ob2}^e(n(0),n(1))$ to be the social benefit per time unit in Case~2, or $N-1\leq n(1) \leq n(0)$. Then,
\begin{equation}\label{bb78}
\begin{split}
U_{ob2}^e  (n(0),n(1)) = \lambda R(1 - {\pi _{0,n(0) + 1}} - \sum\limits_{n = n(1) + 1}^{n(0) + 1} {{\pi _{1,n}}} ) \\
 - C(\sum\limits_{n = 0}^{n(0) + 1} {n{\pi _{0,n}} + \sum\limits_{n = 0}^{n(1) + 1} {n{\pi _{1,n}}} }  + \sum\limits_{n = 1}^{n(0) + 1} {n{\pi _{2,n}}} ).
\end{split}
\end{equation}
Obviously, the equilibrium social benefit is $U_{ob2}^e({n_e}(0),{n_e}(1))$.

\begin{figure*}
\centering
\includegraphics[width=1.0\textwidth]{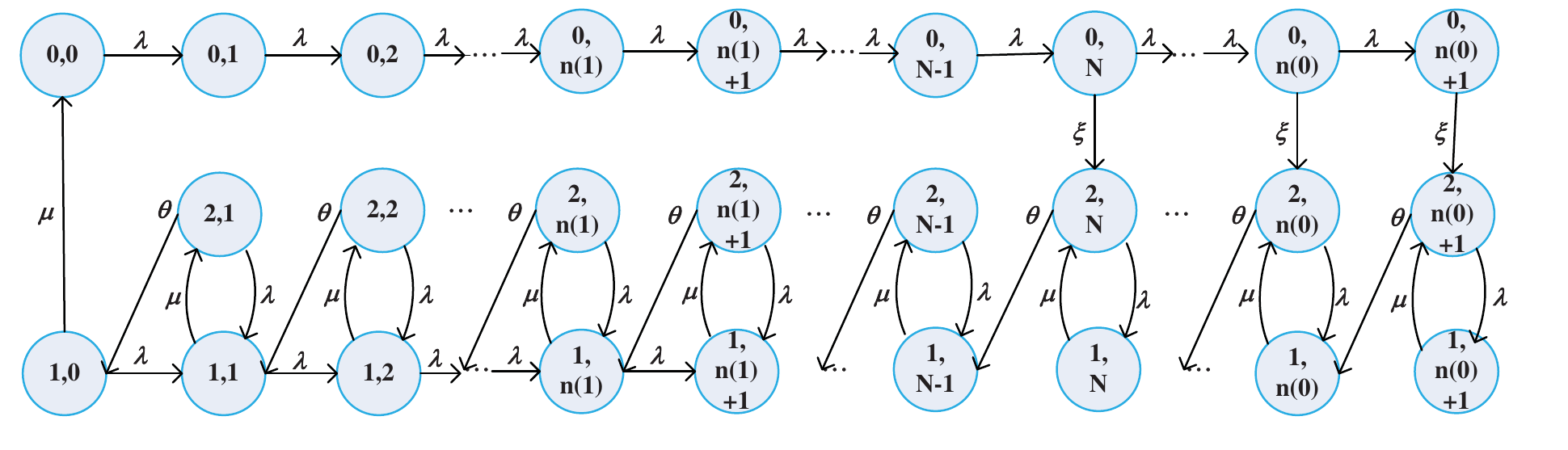}
\caption{Transition rate diagram of {($M(t)$,$I(t)$)} for the observable queues when $n(1) < N - 1 \le n(0)$.}
\label{Fig:4}
\end{figure*}

Case 3: For $n(1) < N - 1 \le n(0)$, the corresponding transition rate diagram is showed in Fig. \ref{Fig:4}, and the state space of $\{ (M(t),I(t))\} $ is given by:
\begin{equation} \label{eqn:case3}
    \Omega _{ob3}^e = \{ (0,n):0 \le n \le n(0) + 1\}  \cup \{ (1,n):0 \le n \le n(0) + 1\}  \cup \{ (2,n):1 \le n \le n(0) + 1\} .
\end{equation}

The stationary distribution for this case (Caese~3) is given in Theorem~\ref{t34}.
\begin{theorem}  \label{t34}
For the fully observable {M/M/1/MV} constant retrial queue with multiple vacations and the $N$-policy, if $n(1) < N - 1 \le n(0)$, then the state space $\Omega _{ob3}^e$ of $\{ (M(t),I(t))\} $ is given by (\ref{eqn:case3}), and the stationary distribution $\{ {\pi _{i,n}}\left| {(i,n) \in \Omega_{ob3}^e} \right.\} $ is given by:
\begin{equation}\label{bb79}
{\pi _{0,n}} =
\begin{cases}
\frac{\mu }{\lambda } \cdot {\pi _{1,0}}, &0 \le n \le N - 1;\\
\frac{\mu }{\lambda } \cdot {\left( {\frac{\lambda }{{\lambda  + \xi }}} \right)^{n - N + 1}} \cdot {\pi _{1,0}}, &N \le n \le n(0);\\
\frac{\mu }{\xi } \cdot {\left( {\frac{\lambda }{{\lambda  + \xi }}} \right)^{n(0) - N + 1}} \cdot {\pi _{1,0}}, &n = n(0) + 1;
\end{cases}
\end{equation}
\begin{equation}\label{bb80}
{\pi _{1,n}} =
\begin{cases}
{A_1} + {A_2} \cdot {F^n}, &0 \le n \le N - 1;\\
{B_1} + {B_2} \cdot {F^n} + {D_1} \cdot {\left( {\frac{\lambda }{{\lambda  + \xi }}} \right)^n}, &N \le n \le n(1);\\
\begin{split}
{B_1} + {B_2} \cdot {F^{n(1) + 1}} + {D_1} \cdot \left( {\frac{{\lambda  - F\xi }}{{\lambda  + \xi }}} \right) \cdot {\left( {\frac{\lambda }{{\lambda  + \xi }}} \right)^{n(1) - 1}}\\
 -(\frac{\xi }{{\lambda  + \xi }} + \frac{\xi }{\theta }) \cdot {\left( {\frac{\lambda }{{\lambda  + \xi }}} \right)^{n(1) - N + 1}} \cdot {\pi _{1,0}},
\end{split}&n = n(1) + 1;\\
(1 - F){B_1} + {D_1}\left( {\frac{\lambda }{{\lambda  + \xi }} - F} \right){\left( {\frac{\lambda }{{\lambda  + \xi }}} \right)^{n(1) - 1}} + \psi (n) \cdot {\pi _{1,0}}, &n(1) + 2 \le n \le n(0)\\
{\pi _{1,n(0)}} - \left( {1 + \frac{\xi }{\theta }} \right){\left( {\frac{\lambda }{{\lambda  + \xi }}} \right)^{n(0) - N + 1}}{\pi _{1,0}}, &n=n(0) + 1;
\end{cases}
\end{equation}
and
\begin{equation}\label{bb81}
{\pi _{2,n}} =
\begin{cases}
\frac{\mu }{{\lambda  + \theta }} \cdot {\pi _{1,n}}, &1 \le n \le N - 1;\\
\frac{\mu }{{\lambda  + \theta }} \cdot {\pi _{1,n}} + \frac{\xi }{{\lambda  + \xi }} \cdot {\pi _{0,n}}, &N \le n \le n(0) + 1;
\end{cases}
\end{equation}
where ${A_i}$ ($i = 1,2)$, ${B_i}$ ($i = 1,2)$, $F$ and ${D_1}$ are given by (\ref{bb19}), (\ref{bb20}), (\ref{bb21}) and (\ref{bb22}), respectively,
\begin{equation}\label{bb82}
\psi (n) = \left( {\frac{\lambda }{{\lambda  + \xi }}} \right.\left. { + \frac{\lambda }{\theta }} \right){\left( {\frac{\lambda }{{\lambda  + \xi }}} \right)^{n - N}} - \left( {1 + \frac{{\lambda  + \xi }}{\theta } + \frac{\xi }{\lambda } + \frac{{\xi (\lambda  + \xi )}}{{\lambda  + \theta }}} \right){\left( {\frac{\lambda }{{\lambda  + \xi }}} \right)^{n(1) + 2 - N}}.
\end{equation}
${\pi _{1,0}}$ can be obtained by the normalization condition $\sum\limits_{(i,n) \in \Omega _{ob3}^e} {{\pi _{i,n}} = 1} $.
\end{theorem}

From Fig. \ref{Fig:4}, we know that Case~3 and Case~2 have the same equilibrium equations. Hence, the stationary distribution of Case~3 is the same as that of Case~2, and therefore we omit the proof of  Theorem \ref{t34}. Based on Fig \ref{Fig:4} and Theorem \ref{t34}, we can obtain that the states, at which customers will balk are $(0,n(0) + 1)$ and $(1,n(1) + 1)$. Denote the social benefit per time unit in Case~3 by $U_{ob3}^e(n(0),n(1))$. We then have
\begin{equation}\label{bb83}
U_{ob3}^e  (n(0),n(1)) = \lambda R(1 - {\pi _{0,n(0) + 1}} - \sum\limits_{n = n(1) + 1}^{n(0) + 1} {{\pi _{1,n}}} )
 - C(\sum\limits_{n = 0}^{n(0) + 1} {n{\pi _{0,n}} + \sum\limits_{n = 0}^{n(1) + 1} {n{\pi _{1,n}}} }  + \sum\limits_{n = 1}^{n(0) + 1} {n{\pi _{2,n}}} ).
\end{equation}
Obviously, the equilibrium social benefit is $U_{ob3}^e({n_e}(0),{n_e}(1))$.
\subsection{Social optimization}
Summarize the above discussion, we define ${U_s}(n(0),n(1))$ as the social benefit per time unit, thus
\begin{equation}\label{bb831}
{U_s}(n(0),n(1)) =
\begin{cases}
U_{ob1}^e(n(0),n(1)) &{\rm if}~~N - 1 \le n(0) \le n(1);\\
U_{ob2}^e(n(0),n(1)) &{\rm if}~~N - 1 \le n(1) \le n(0);\\
U_{ob3}^e(n(0),n(1)) &{\rm if}~~n(1) < N - 1 \le n(0).
\end{cases}
\end{equation}
We define ${U_s}({n^*}(0),{n^*}(1))$ as the socially optimal social welfare. Obviously, it's easy for us to get ${U_s}({n^*}(0),{n^*}(1)) = max\{ {U_s}(n(0),n(1))\} $.

\section{The unobservable case}
\label{sec:4}

\begin{figure*}
\centering
\includegraphics[width=0.8\textwidth]{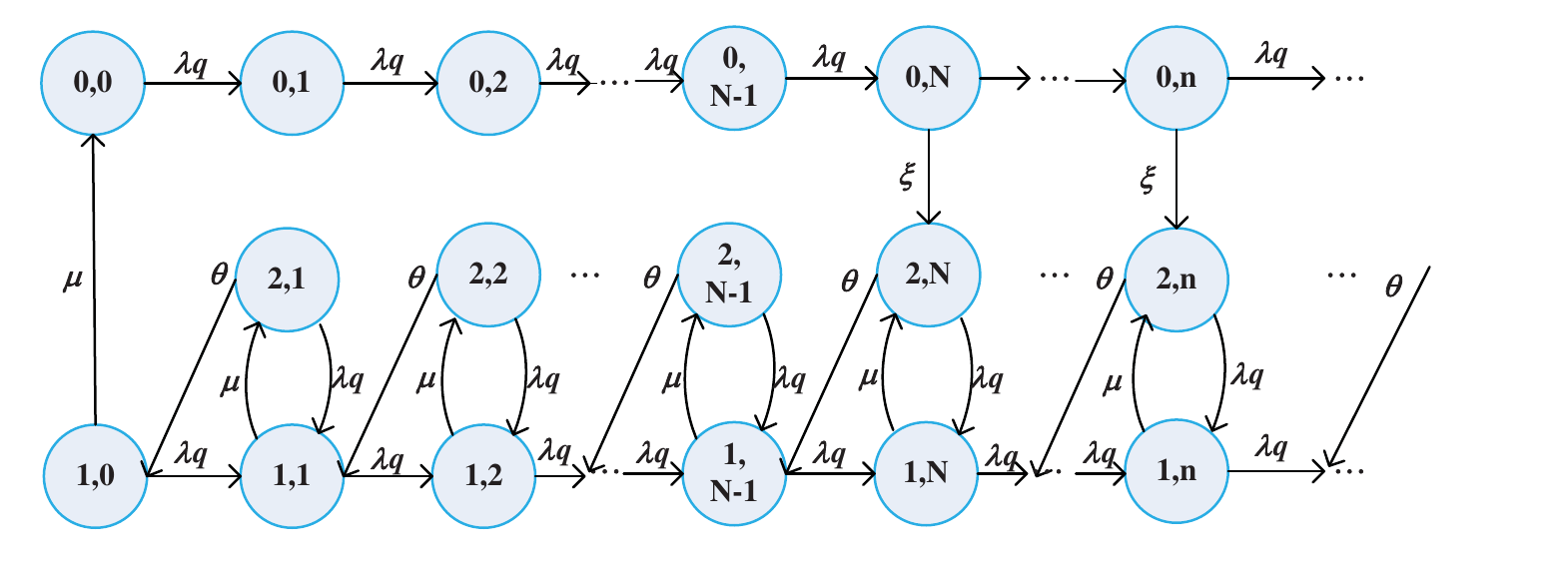}
\caption{Transition rate diagram for the unobservable case.}
\label{Fig:5}
\end{figure*}

In the observable case, we assume that the arriving customers can observe all information about $M(t)$ and $I(t)$. Now, in the unobservable case, we assume that the arriving customers cannot observe any information about $M(t)$ or $I(t)$. In this case, the customers join the system with probability $q$ ($0 \le q \le 1$). So, the effective arrival rate is $\overline \lambda=\lambda q$, the equilibrium mixed strategy of the customers is denoted by the equilibrium arrival rate ${\overline \lambda  _e}=\lambda {q_e}$ (${q_e}$ is equilibrium joining probability of the customers), and the socially optimal mixed strategy is denoted by the optimal arrival rate ${\overline \lambda  ^*}=\lambda {q^*}$, where ${q^*}$ is optimal joining probability of the customers. The corresponding transition rate diagram is showed in Fig. \ref{Fig:5}.

\subsection{Equilibrium}

In the unobservable case, the arriving customers can neither observe the state of the server $M(t)$ \red{nor} the number of other customers $I(t)$ in the orbit. In order to obtain the equilibrium arrival rate ${\overline \lambda  _e}$ in this case, the stationary distribution needs to be determined first. From Fig. \ref{Fig:5}, it is easy to know that the process $\{M(t),I(t)\}$ is a quasi-birth-and-death (QBD) process with
$${\Omega _{un}} = \{ (0,n):n \ge 0\}  \cup \{ (1,n):n \ge 0\}  \cup \{ (2,n):n \ge 1\} .$$
If $\overline \rho   = \overline \lambda  /\mu  < 1$, $(M,I)$ is defined as the stationary limit of the process $\{M(t),I(t)\}$. The stationary distribution is denoted by:
\[
\pi  = (\pi _{0,0},\pi_{1,0}, \pi_1, \pi_2, \ldots, \pi_n, \ldots ),
\]
where
\[
    \pi_n = ({\pi _{0,n}},{\pi _{1,n}},{\pi _{2,n}}),~~n \ge 1,
\]
with the definition of
\[
    {\pi _{i,n}} = P\left\{ {M = i,I = n} \right\} = \mathop {\lim }\limits_{t \to \infty } P\left\{ {M(t) = i,I(t) = n} \right\}, \quad (i,n) \in {\Omega _{un}} \mbox{ for } i = 0,1,2.
\]
Let $Q$ be defined as the infinitesimal generator of the process. Then, the stationary probability vector $\pi$ can be solved through the equations $\pi Q = 0$:
\begin{align}\label{cc1}
    \overline \lambda  {\pi _{0,0}} &= \mu {\pi _{1,0}};\\
\label{cc2}
    \overline \lambda  {\pi _{0,n}} &= \overline \lambda  {\pi _{0,n - 1}}, \qquad  1 \le n \le N - 1;\\
\label{cc3}
    (\overline \lambda   + \xi ){\pi _{0,n}} &= \overline \lambda  {\pi _{0,n - 1}}, \qquad  n \ge N;\\
\label{cc4}
    (\overline \lambda   + \mu ){\pi _{1,0}} &= \theta {\pi _{2,1}};\\
\label{cc5}
    (\overline \lambda   + \mu ){\pi _{1,n}} &= \overline \lambda  {\pi _{1,n - 1}} + \overline \lambda  {\pi _{2,n}} + \theta {\pi _{2,n + 1}}, \qquad  n \ge 1;\\
\label{cc6}
    (\overline \lambda   + \theta ){\pi _{2,n}} &= \mu {\pi _{1,n}}, \qquad  1 \le n \le N - 1;\\
\label{cc7}
    (\overline \lambda   + \theta ){\pi _{2,n}} &= \xi {\pi _{0,n}} + \mu {\pi _{1,n}}, \qquad  n \ge N.
\end{align}

From Fig. \ref{Fig:5} and ordering the states in the state space ${\Omega _{un}}$ lexicographically,  we can write the infinite generator for the Markov process $\{ (M(t),I(t))\} $ as a block-partitioned form as follows:
\begin{equation}\label{cc8}
 Q = \left( {\begin{array}{*{20}{c}}
{{A_{0,0}}}&{{A_{0,1}}}&{}&{}&{}&{}&{}&{}&{}&{}\\
{{A_{1,0}}}&{{A_1}}&{{A_0}}&{}&{}&{}&{}&{}&{}&{}\\
{}&{{A_2}}&{{A_1}}&{{A_0}}&{}&{}&{}&{}&{}&{}\\
{}&{}&{{A_2}}&{{A_1}}&{{A_0}}&{}&{}&{}&{}&{}\\
{}&{}&{}& \ddots & \ddots & \ddots &{}&{}&{}&{}\\
{}&{}&{}&{}&{{A_2}}&{{A_1}}&{{A_0}}&{}&{}&{}\\
{}&{}&{}&{}&{}&{{A_2}}&{{B_1}}&{{A_0}}&{}&{}\\
{}&{}&{}&{}&{}&{}&{{A_2}}&{{B_1}}&{{A_0}}&{}\\
{}&{}&{}&{}&{}&{}&{}& \ddots & \ddots & \ddots
\end{array}} \right),
\end{equation}
where
\begin{align*}
    {A_{0,0}} &= \left( {\begin{array}{*{20}{c}}
{ - \overline \lambda  }&0\\
\mu &{ - (\overline \lambda   + \mu )}
\end{array}} \right), \qquad
{A_{0,1}} = \left( {\begin{array}{*{20}{c}}
{\overline \lambda  }&0&0\\
0&{\overline \lambda  }&0
\end{array}} \right), \qquad
{A_{1,0}} = \left( {\begin{array}{*{20}{c}}
0&0\\
0&0\\
0&\theta
\end{array}} \right), \\
    {A_0} &= \left( {\begin{array}{*{20}{c}}
{\overline \lambda  }&0&0\\
0&{\overline \lambda  }&0\\
0&0&0
\end{array}} \right), \qquad
{A_1} = \left( {\begin{array}{*{20}{c}}
{ - \bar \lambda }&0&0\\
0&{ - (\bar \lambda  + \mu )}&\mu \\
0&{\bar \lambda }&{ - (\bar \lambda  + \theta )}
\end{array}} \right), \qquad
{A_2} = \left( {\begin{array}{*{20}{c}}
0&0&0\\
0&0&0 \\
0&\theta&0
\end{array}} \right), \\
    {B_1} &= \left( {\begin{array}{*{20}{c}}
{ - (\bar \lambda  + \xi )}&0&\xi \\
0&{ - (\bar \lambda  + \mu )}&\mu \\
0&{\bar \lambda }&{ - (\bar \lambda  + \theta )}
\end{array}} \right).
\end{align*}

\begin{theorem}  \label{t411}
For the unobservable case, the {M/M/1/MV} constant retrial queue with multiple vacations and the $N$-policy, given the arrival rate ${\overline \lambda  }$, the stationary distribution $\{ {\pi _{i,n}}\left| {(i,n) \in \Omega _{un}} \right.\} $ is given by
\begin{equation}\label{cc9}
{\pi _{0,n}} =
\begin{cases}
\frac{\mu }{{\bar \lambda }}{\pi _{1,0}}, &0 \le n \le N - 1;\\
\frac{\mu }{{\bar \lambda  + \xi }}{\pi _{1,0}}, &n = N;\\
{r_{11}}{\pi _{0,N}}, &n \geq N;\\
\end{cases}
\end{equation}
\begin{equation}\label{cc10}
{\pi _{1,n}} =
\begin{cases}
\overline {{A_1}}  + \overline {{A_2}}  \cdot {\overline F ^n},0 \le n \le N - 1, &0 \le n \le N - 1;\\
\frac{{\overline \rho  (\overline \lambda   + \theta )}}{\theta }(\overline {{A_1}}  + \overline {{A_2}} \cdot {\overline F ^{N - 1}}) + \frac{{\overline \lambda  (\overline \lambda   + \theta )}}{{\theta (\overline \lambda   + \xi )}}{\pi _{1,0}} + \frac{{\overline \lambda  \xi }}{{\theta (\overline \lambda   + \xi )}}{\pi _{1,0}}, &n = N;\\
{r_{12}}{\pi _{0,N}} + {r_{22}}{\pi _{1,N}}, &n \geq N;
\end{cases}
\end{equation}
and
\begin{equation}\label{cc11}
{\pi _{2,n}} =
\begin{cases}
\frac{\mu }{{\bar \lambda  + \theta }}{\pi _{1,n}}, &1\le n \le N - 1;\\
\frac{\xi }{{\bar \lambda  + \theta }}{\pi _{0,N}} + \frac{\mu }{{\bar \lambda  + \theta }}{\pi _{1,N}}, &n = N;\\
{r_{13}}{\pi _{0,N}} + {r_{23}}{\pi _{1,N}},  &n \geq N;
\end{cases}
\end{equation}
where $\overline F  = \frac{{\overline \lambda  (\overline \lambda   + \theta )}}{{\theta \mu }}$,
\begin{equation}\label{cc12}
\left\{ \begin{array}{l}
\overline {{A_1}}  = \frac{{\mu \overline F }}{{\overline \lambda  (1 - \overline F )}}{\pi _{1,0}}; \\
\overline {{A_2}}  = (1 + \frac{{\mu \overline F }}{{\overline \lambda  (\overline F  - 1)}}){\pi _{1,0}};
\end{array} \right.
\end{equation}
\begin{equation}\label{cc13}
 \left\{ \begin{array}{l}
{r_{11}} = {\left( {\frac{{\bar \lambda }}{{\bar \lambda  + \xi }}} \right)^n};\\
{r_{12}} = \frac{{\bar \lambda {{\left( {\frac{{\bar \lambda (\bar \lambda  + \theta )}}{{\theta \mu }}} \right)}^{n - N}}(\bar \lambda  + \theta  + \xi )}}{{\bar \lambda \theta  + {{\bar \lambda }^2} - \theta \mu  + \theta \xi  + \bar \lambda \xi }} - \frac{{{{\left( {\frac{{\bar \lambda }}{{\bar \lambda  + \xi }}} \right)}^{n - N}}({{\bar \lambda }^2} + \bar \lambda \theta  + \bar \lambda \xi )}}{{\bar \lambda \theta  + {{\bar \lambda }^2} - \theta \mu  + \theta \xi  + \bar \lambda \xi }};\\
{r_{13}} = \frac{{\bar \lambda {{\left( {\frac{{\bar \lambda (\bar \lambda  + \theta )}}{{\theta \mu }}} \right)}^{n - N}}\mu (\bar \lambda  + \theta  + \xi )}}{{(\bar \lambda  + \theta )(\bar \lambda \theta  + {{\bar \lambda }^2} - \theta \mu  + \theta \xi  + \bar \lambda \xi )}} - \frac{{{{\left( {\frac{{\bar \lambda }}{{\bar \lambda  + \xi }}} \right)}^{n - N}}(\bar \lambda \mu  - \bar \lambda \xi  + \mu \xi  - {\xi ^2})}}{{\bar \lambda \theta  + {{\bar \lambda }^2} - \theta \mu  + \theta \xi  + \bar \lambda \xi }};\\
{r_{22}} = {\left( {\frac{{\bar \lambda (\bar \lambda  + \theta )}}{{\theta \mu }}} \right)^{n - N}};\\
{r_{23}} = \frac{{\mu {{\left( {\frac{{\bar \lambda (\bar \lambda  + \theta )}}{{\theta \mu }}} \right)}^{n - N}}}}{{\bar \lambda  + \theta }};
\end{array} \right.
\end{equation}
and
\begin{equation}\label{cc14}
{\pi _{1,0}} = \frac{{\overline \lambda  \xi (\theta \mu  - {{\overline \lambda  }^2} - \overline \lambda  \theta )}}{{{\mu ^2}(\overline \lambda   + \theta )(\overline \lambda   + N\xi )}}.
\end{equation}
\end{theorem}

{\bf Proof}~~In order to obtain the stationary distribution of the system, we first need to obtain the rate matrix $R$, which is the minimum non-negative solution of the following matrix quadratic equation:
\begin{equation}\label{cc15}
{R^2}{A_2} + R{B_1} + {A_0} = 0.
\end{equation}
By detailed calculations, we get the minimum non-negative solution of $R$ as follows:
\begin{equation}\label{cc16}
R = \left( {\begin{array}{*{20}{c}}
{\frac{{\overline \lambda  }}{{\bar \lambda  + \xi }}}&{\frac{{{{\bar \lambda }^2}(\bar \lambda  + \theta  + \xi )}}{{\theta \mu (\bar \lambda  + \xi )}}}&{\frac{{\bar \lambda }}{\theta }}\\
0&{\frac{{\bar \lambda (\bar \lambda  + \theta )}}{{\theta \mu }}}&{\frac{{\bar \lambda }}{\theta }}\\
0&0&0
\end{array}} \right).
\end{equation}
Using the matrix-geometric solution (see \cite{neuts1981matrix}), we have:
\begin{equation}\label{cc17}
{\pi _n} = {\pi _N}{R^{n - N}},n \ge N,
\end{equation}
and $({\pi _{0,0}},{\pi _{1,0}},{\pi _1},{\pi _2} \ldots {\pi _{N - 1}},{\pi _N})$ satisfies the following equation:
\begin{equation}\label{cc18}
({\pi _{0,0}},{\pi _{1,0}},{\pi _1},{\pi _2} \ldots {\pi _{N - 1}},{\pi _N})B[R] = 0,
\end{equation}
where
\begin{equation}\label{cc19}
B[R] = \left( {\begin{array}{*{20}{c}}
{{A_{0,0}}}&{{A_{0,1}}}&{}&{}&{}&{}&{}\\
{{A_{1,0}}}&{{A_1}}&{{A_0}}&{}&{}&{}&{}\\
{}&{{A_2}}&{{A_1}}&{{A_0}}&{}&{}&{}\\
{}&{}&{{A_2}}&{{A_1}}&{{A_0}}&{}&{}\\
{}&{}&{}& \ddots & \ddots & \ddots &{}\\
{}&{}&{}&{}&{{A_2}}&{{A_1}}&{{A_0}}\\
{}&{}&{}&{}&{}&{{A_2}}&{R{A_2} + {B_1}}
\end{array}} \right).
\end{equation}
Substituting (\ref{cc19}) into (\ref{cc18}), we can get:
\begin{equation}\label{cc20}
\left\{ \begin{array}{l}
\overline \lambda  {\pi _{0,0}} = \mu {\pi _{1,0}};\\
\bar \lambda {\pi _{0,n}} = \bar \lambda {\pi _{0,n - 1}}; \qquad 1 \le n \le N - 1;\\
(\bar \lambda  + \xi ){\pi _{0,N}} = \bar \lambda {\pi _{0,N - 1}};\\
(\bar \lambda  + \mu ){\pi _{1,0}} = \theta {\pi _{2,1}};\\
(\bar \lambda  + \mu ){\pi _{1,n}} = \bar \lambda {\pi _{1,n - 1}} + \bar \lambda {\pi _{2,n}} + \theta {\pi _{2,n + 1}}, \qquad 1 \le n \le N - 1;\\
\mu {\pi _{1,N}} = \bar \lambda {\pi _{1,N - 1}} + \bar \lambda {\pi _{0,N}} + \bar \lambda {\pi _{2,N}};\\
(\bar \lambda  + \theta ){\pi _{2,n}} = \mu {\pi _{1,n}}, \qquad 1 \le n \le N - 1;\\
(\bar \lambda  + \theta ){\pi _{2,N}} = \xi {\pi _{0,N}} + \mu {\pi _{1,N}}.
\end{array} \right.
\end{equation}
By calculating (\ref{cc20}), we can get
\begin{equation}\label{cc21}
{\pi _{0,n}} =
\begin{cases}
\frac{\mu }{{\bar \lambda }}{\pi _{1,0}}, &0 \le n \le N - 1;\\
\frac{\mu }{{\bar \lambda  + \xi }}{\pi _{1,0}}, &n = N;
\end{cases}
\end{equation}
\begin{equation}\label{cc22}
{\pi _{1,n}} =
\begin{cases}
\overline {{A_1}}  + \overline {{A_2}}\cdot {\overline F ^n}, &0 \le n \le N - 1;\\
\frac{{\overline \rho  (\overline \lambda   + \theta )}}{\theta }(\overline {{A_1}}  + \overline {{A_2}} \cdot {\overline F ^{N - 1}}) + \frac{{\overline \lambda  (\overline \lambda   + \theta )}}{{\theta (\overline \lambda   + \xi )}}{\pi _{1,0}} + \frac{{\overline \lambda  \xi }}{{\theta (\overline \lambda   + \xi )}}{\pi _{1,0}}, &n = N;
\end{cases}
\end{equation}
and
\begin{equation}\label{cc23}
{\pi _{2,n}} =
\begin{cases}
\frac{\mu }{{\bar \lambda  + \theta }}{\pi _{1,n}}, &0 \le n \le N - 1;\\
\frac{\xi }{{\bar \lambda  + \theta }}{\pi _{0,N}} + \frac{\mu }{{\bar \lambda  + \theta }}{\pi _{1,N}}, &n = N;
\end{cases}
\end{equation}
 where $\overline F  = \frac{{\overline \lambda  (\overline \lambda   + \theta )}}{{\theta \mu }}$,
\begin{equation}\label{cc24}
\left\{ \begin{array}{l}
\overline {{A_1}}  = \frac{{\mu \overline F }}{{\overline \lambda  (1 - \overline F )}}{\pi _{1,0}},\\
\overline {{A_2}}  = (1 + \frac{{\mu \overline F }}{{\overline \lambda  (\overline F  - 1)}}){\pi _{1,0}}.
\end{array} \right.
\end{equation}

From (\ref{cc21})\red{--}(\ref{cc23}), ${\pi _N} = ({\pi _{0,N}},{\pi _{1N}},{\pi _{2,N}})$ can be obtained. By (\ref{cc16}), ${R^{n - N}}$ can be obtained as follows:
\begin{equation}\label{cc25}
{R^{n - N}}=\left( {\begin{array}{*{20}{c}}
{{r_{11}}}&{{r_{12}}}&{{r_{13}}}\\
0&{{r_{22}}}&{{r_{23}}}\\
0&0&0
\end{array}} \right),
\end{equation}
where ${{r_{11}}}, {{r_{12}}}, {{r_{13}}}, {{r_{22}}}$ and ${{r_{23}}}$ see (\ref{cc13}). Now, from (\ref{cc17}), we can get:
\begin{equation}\label{cc26}
\begin{cases}
{\pi _{0,n}} = {r_{11}}{\pi _{0,N}}, &n \ge N;\\
{\pi _{1,n}} = {r_{12}}{\pi _{0,N}} + {r_{22}}{\pi _{1,N}}, &n \ge N;\\
{\pi _{2,n}} = {r_{13}}{\pi _{0,N}} + {r_{23}}{\pi _{1,N}},  &n \ge N.
\end{cases}
\end{equation}
Hence, considering the above discussion, (\ref{cc9})--(\ref{cc11}) can be obtained, and ${\pi _{1,0}}$ can be calculated by the following normalization condition:
\begin{equation}\label{cc27}
{\pi _{0,0}} + {\pi _{1,0}} + \sum\limits_{n = 1}^{N - 1} {{\pi _n}e}  + {\pi _N}{(I - R)^{ - 1}}e = 1,
\end{equation}
where the expression for ${\pi_{1,0}}$ is given in (\ref{cc14}). \pend

In Theorem \ref{t411}, the stationary distribution under unobservable case was obtained by using the matric-analytic method, based on which we can get the mean queue length $E[L(\bar \lambda )]$ for the unobservable case, given by:
\begin{align}\label{cc28}
  \nonumber E[L(\bar \lambda )] & = \sum\limits_{n = 1}^\infty  {n({\pi _{0,n}} + {\pi _{1,n}} + {\pi _{2,n}})} \\
   & =\frac{{\bar \lambda }}{\xi } + \frac{{\bar \lambda }}{{\theta (\mu  - \bar \lambda )}} + \frac{{\theta \xi }}{{(\bar \lambda  + \theta )(\bar \lambda  + \xi )}} + \frac{{\theta (\mu  - \xi )}}{{{\mu ^2}(\bar \lambda  + \xi )}} + \frac{{N(N - 1)\xi }}{{{\mu ^2}(\bar \lambda  + N\xi )}}.
\end{align}
By using Little's law for the whole system, we can get mean sojourn time $E[W(\bar \lambda )]$:
\begin{equation}\label{cc29}
E[W(\overline \lambda  )] = \frac{{E[L(\overline \lambda  )]}}{{\overline \lambda  }} = \frac{1}{\xi } + \frac{1}{{\theta (\mu  - \bar \lambda )}} + \frac{{\theta \xi }}{{\bar \lambda (\bar \lambda  + \theta )(\bar \lambda  + \xi )}} + \frac{{\theta (\mu  - \xi )}}{{\bar \lambda {\mu ^2}(\bar \lambda  + \xi )}} + \frac{{N(N - 1)\xi }}{{\bar \lambda {\mu ^2}(\bar \lambda  + N\xi )}}.
\end{equation}
We obtain the equilibrium arrival rate and socially optimal arrival rate by the following theorem.

\begin{theorem}  \label{t412}
For the unobservable {M/M/1/MV} constant retrial queue with multiple vacations and the $N$-policy, we have the following conclusions on the equilibrium arrival rate: \\
(i) It has no positive equilibrium arrival rate, if $R < CE[W(\bar \lambda )]$;\\
(ii) It has one positive equilibrium arrival rate $\overline {{\lambda _e}}  = \overline {{\lambda _1}} $ (if and only if $\overline {{\lambda _1}}  \le \lambda $), if $R = CE[W(\bar \lambda )]$, where $\overline {{\lambda _1}}$ is the unique positive solution of $E[W'(\overline \lambda  )] = 0$; \\
(iii)
\begin{equation}\label{cc30}
\overline {{\lambda _e}}
\begin{cases}
\in \{ \overline {{\lambda _2}} ,\overline {{\lambda _3}} \}, & {\rm if}~~ R>CE[W(\overline {{\lambda _1}} )] ~~{\rm and}~~ \overline {{\lambda _3}}  \le \lambda;\\
\in \{ \overline {{\lambda _2}} ,\lambda \}, & {\rm if}~~ R>CE[W(\overline {{\lambda _1}} )] ~~{\rm and}~~ \overline {{\lambda _2}}  < \lambda  < \overline {{\lambda _3}};\\
= \lambda, & {\rm if}~~ R>CE[W(\overline {{\lambda _1}} )] ~~{\rm and}~~ \lambda  = \overline {{\lambda _2}} ;\\
\rm{no~~ positive ~~rquilibrium ~~rate}, & {\rm if}~~ R>CE[W(\overline {{\lambda _1}} )] ~~{\rm and}~~\overline {{\lambda _2}}  > \lambda;
\end{cases}
\end{equation}
where $\overline {{\lambda _2}}$, $\overline {{\lambda _3}}$ (with $0 \le \overline {{\lambda _2}}  \le \overline {{\lambda _3}} $) are the positive solutions of $R - CE[W(\overline \lambda  )] = 0$.
\end{theorem}

{\bf Proof}~~From (\ref{cc29}), we can obtain the expected net benefit ${U_e}(\overline \lambda  )$ of the marked customer:
\begin{align} \label{cc31}
    \nonumber {U_e}(\overline \lambda  ) &= R - CE[W(\overline \lambda  )]\\
    & = R - C\left( {\frac{1}{\xi } + \frac{1}{{\theta (\mu  - \bar \lambda )}} + \frac{{\theta \xi }}{{\bar \lambda (\bar \lambda  + \theta )(\bar \lambda  + \xi )}} + \frac{{\theta (\mu  - \xi )}}{{\bar \lambda {\mu ^2}(\bar \lambda  + \xi )}} + \frac{{N(N - 1)\xi }}{{\bar \lambda {\mu ^2}(\bar \lambda  + N\xi )}}} \right).
\end{align}
Since the second-order derivative of $E[W(\overline \lambda  )]$ in ${\overline \lambda  }$ is given by:
\begin{align}\label{cc32}
  \nonumber E[W''(\overline \lambda  )] &  = \frac{2}{{\theta   {{(\mu  - \bar \lambda )}^3}}} + 2\theta \xi \left( {\frac{1}{{{{\bar \lambda }^3}(\bar \lambda  + \theta )(\bar \lambda  + \xi )}} + \frac{1}{{\bar \lambda (\bar \lambda  + \theta ){{(\bar \lambda  + \xi )}^3}}}} \right. \\
  \nonumber &   + \left. {\frac{1}{{\bar \lambda {{(\bar \lambda  + \theta )}^2}{{(\bar \lambda  + \xi )}^2}}} + \frac{1}{{\bar \lambda {{(\bar \lambda  + \theta )}^3}(\bar \lambda  + \xi )}} + \frac{1}{{{{\bar \lambda }^2}(\bar \lambda  + \theta ){{(\bar \lambda  + \xi )}^2}}} + \frac{1}{{{{\bar \lambda }^2}{{(\bar \lambda  + \theta )}^2}(\bar \lambda  + \xi )}}} \right) \\
  \nonumber &  + 2\theta (\mu  - \xi )\left( {\frac{1}{{\bar \lambda {\mu ^2}{{(\bar \lambda  + \xi )}^3}}} + \frac{1}{{{{\bar \lambda }^2}{\mu ^2}{{(\bar \lambda  + \xi )}^2}}} + \frac{1}{{{{\bar \lambda }^3}{\mu ^3}(\bar \lambda  + \xi )}}} \right) \\
   &  + 2N(N - 1)\xi \left( {\frac{1}{{\bar \lambda {\mu ^2}(\bar \lambda  + N\xi )}} + \frac{1}{{{{\bar \lambda }^2}{\mu ^2}{{(\bar \lambda  + N\xi )}^2}}} + \frac{1}{{{{\bar \lambda }^3}{\mu ^2}(\bar \lambda  + N\xi )}}} \right).
\end{align}
If $\overline \rho   = \frac{{\overline \lambda  }}{\mu } < 1,$ then $E[W''(\overline \lambda  )]$ is positive. In this case, $E[W(\overline \lambda  )]$ is strictly convex of $\overline \lambda$. If we denote the positive solution of equation $E[W'(\overline \lambda  )] = 0$ by $\overline {{\lambda _1}} $, and the positive solutions of ${U_e}(\overline \lambda  ) = R - CE[W(\overline \lambda  )] = 0$ by $\overline {{\lambda _2}}$ and $\overline {{\lambda _3}} $ ($0\leq \overline {{\lambda _2}} \leq \overline {{\lambda _3}} $), then we have the following results:

(1) When $R < CE[W(\overline {{\lambda _1}} )]$, i.e., ${U_e}(\overline {{\lambda _1}} ) < 0$, ${U_e}(\overline \lambda  )$ is negative for every $\overline \lambda$. Therefore, it has no positive equilibrium arrival rate, which leads to (i).

(2) When $R = CE[W(\overline {{\lambda _1}} )]$, i.e., ${U_e}(\overline {{\lambda _1}} ) = 0$, ${U_e}(\overline \lambda  )$ is negative for every $\overline \lambda   \ne \overline {{\lambda _1}} $. Therefor, it has one positive equilibrium arrival rate $\overline {{\lambda _e}}  = \overline {{\lambda _1}} $ (if and only if $\overline {{\lambda _1}}  \le \lambda $), which leads to (ii).

(3) When $R > CE[W(\overline {{\lambda _1}} )]$, i.e., ${U_e}(\overline {{\lambda _1}} ) > 0$. If $\overline {{\lambda _3}}  \le \lambda $, it has two positive equilibrium arrival rates $\overline {{\lambda _e}}  = \{ \overline {{\lambda _2}} ,\overline {{\lambda _3}} \} $, which leads to the first part of (\ref{cc30}). If $\overline {{\lambda _2}}  < \lambda  < \overline {{\lambda _3}}  $, it has two positive equilibrium arrival rates $\overline {{\lambda _e}}  = \{ \overline {{\lambda _2}} ,\lambda \}  $, which leads to the second part of (\ref{cc30}). If $\lambda  = \overline {{\lambda _2}} $, it has a unique positive equilibrium arrival rate $\overline {{\lambda _e}} =\overline {{\lambda _2}} $, which leads to the third part of (\ref{cc30}). If $\overline {{\lambda _2}}  > \lambda  $, obviously, it has no positive equilibrium arrival rate, which leads to the fourth part of (\ref{cc30}).
\pend

\subsection{Social optimization}

In this section, we discuss the optimal balking behavior of customers in the unobservable case in terms of the following Theorem.
\begin{theorem}
For the unobservable {M/M/1/MV} constant retrial queue with multiple vacations and the $N$-policy, the socially optimal mixed strategy is given by
\begin{equation}\label{cc33}
{\overline \lambda  ^*}=
\begin{cases}
\overline \lambda  _1^*, & {\rm if}~~\overline \lambda  _1^* \le \lambda  ;\\
\lambda , &{\rm if}~~ \overline \lambda  _1^* > \lambda  ;
\end{cases}
\end{equation}
where $ \overline \lambda  _1^* $ ($\overline \lambda  _1^* >0$) is the solution of ${U_S}^\prime (\overline \lambda  ) = 0$, and ${U_s}(\overline \lambda  ) = \overline \lambda  R - CE[L(\overline \lambda  )]$.
\end{theorem}

{\bf Proof}~~From (\ref{cc28}), we can get the social welfare per time unite ${U_s}(\overline \lambda  )$ as follows:
\begin{align} \label{cc34}
    \nonumber {U_s}(\overline \lambda  ) &=\overline \lambda  R - CE[L(\overline \lambda  )]\\
    & = \overline \lambda R - C\left( \frac{{\bar \lambda }}{\xi } + \frac{{\bar \lambda }}{{\theta (\mu  - \bar \lambda )}} + \frac{{\theta \xi }}{{(\bar \lambda  + \theta )(\bar \lambda  + \xi )}} + \frac{{\theta (\mu  - \xi )}}{{{\mu ^2}(\bar \lambda  + \xi )}} + \frac{{N(N - 1)\xi }}{{{\mu ^2}(\bar \lambda  + N\xi )}} \right).
\end{align}
The second-order derivative of ${U_s}(\overline \lambda  )$ in ${\overline \lambda  }$ is given by:
\begin{align} \label{cc35}
    \nonumber {U_s}^{\prime \prime }(\overline \lambda  ) & =  - \frac{{2C}}{{\theta {{(\mu  - \overline \lambda  )}^2}}} - \frac{{2C\overline \lambda  }}{{\theta {{(\mu  - \overline \lambda  )}^3}}} - \frac{{2C\theta \xi }}{{(\overline \lambda   + \theta ){{(\overline \lambda   + \xi )}^3}}} - \frac{{2C\theta \xi }}{{{{(\overline \lambda   + \theta )}^2}{{(\overline \lambda   + \xi )}^2}}}\\
    &  - \frac{{2C\theta \xi }}{{{{(\overline \lambda   + \theta )}^3}(\overline \lambda   + \xi )}} - \frac{{2C\theta (\mu  - \xi )}}{{{\mu ^2}{{(\overline \lambda   + \xi )}^3}}} - \frac{{2CN(N - 1)\xi }}{{{\mu ^2}{{(\overline \lambda   + N\xi )}^3}}}.
\end{align}
If $\overline \rho   = \frac{{\overline \lambda  }}{\mu } < 1$, then ${U_s}^{\prime \prime }(\overline \lambda  ) < 0$. So, ${U_s}(\overline \lambda  )$ is strictly concave of $\overline \lambda$. If we denote the unique positive solution of equation ${U_s}^\prime (\overline \lambda  ) = 0$ by $\overline \lambda  _1^* $ ($\overline \lambda  _1^* > 0$), then we have the following results:

(1) When $\overline \lambda  _1^* \le \lambda  $, obviously, the socially optimal mixed strategy ${\overline \lambda  ^*}$ of customers is unique, which is $\overline \lambda  _1^*$.

(2) When $\overline \lambda  _1^* > \lambda $, the socially optimal mixed strategy ${\overline \lambda  ^*}$ of customers is unique, which is $\lambda$. \pend

\section{Numerical results}
\label{sec:5}

In this section, we explore the previous theoretical results through numerical experiments. One should note that, due to the complexity of equations (\ref{bb831}) and (\ref{cc34}), explicit expressions for the equilibrium balking threshold ${n_e}(i)$, $(i=0,1)$, of customers, the socially optimal balking threshold ${n^*}(i)$, $(i=0,1)$, and the optimal social welfare are not available in general. Hence, we use Particle Swarm Optimization (PSO) algorithm to numerically solve complex analytic characteristics in this section. The numerical optimal solution $({n^*}(1),{n^*}(2))$ of $\mathop {\max }\limits_{(n(0),n(1))} {U_s}(n(0),n(1))$ and optimal social welfare ${U_s}({n^*}(0),{n^*}(1))$ and ${U_s}({\overline \lambda  ^*})$ can be obtained by PSO algorithm. PSO algorithm was  introduced by Kennedy and Eberhart~\cite{eberhart1995new} to solve continuous nonlinear optimization problems, and it has been widely used to solve global optimal solutions, since  it does not require many constraints and objective functions. The key procedure of applying PSO algorithm to find the optimal solution (searching for socially optimal balking threshold ${n^*}(i)$ $(i=0,1)$) is illustrated in Algorithm \ref{t1}, where the velocity ${V_{id}}$ and position ${X_{id}}$ are generally provided by:
\begin{equation}\label{c1}
 {V_{id}} = \omega *{V_{id}} + {c_1}*rand()*(pBes{t_i} - {X_{id}}) + {c_2}*rand()*(gBes{t_i} - {X_{id}}),
\end{equation}
and
\begin{equation}\label{c2}
{X_{id}} = {X_{id}} + {V_{id}},
\end{equation}
where $i$ is the number of particles, $d$ is the dimension, $rand()$ is a random number in $(0,1)$, $c_{1}$ and $c_{2}$ are the learning factor and $\omega$ is the inertia factor.
 \begin{algorithm}[htb]
    \caption{Searching for socially optimal balking threshold ${n^*}(0)$ and ${n^*}(1)$}
    \label{t1}
    \begin{algorithmic}[1] 
    \REQUIRE ~~$R$, $C$, $\lambda$, $\mu$, $\xi$, $\theta$, $N$;\\ 
    \ENSURE ~~${n^*}(0)$ , ${n^*}(1)$;\\ 
    \STATE \textbf{for} each particle $i$\\   
    \STATE ~~~~Initializing velocity ${V_{id}}$ and position ${X_{id}}$ for each particle $i$\\
    \STATE ~~~~Evaluating particle $i$ and setting $pBest_{i}={X_{id}}$\\
    \STATE \textbf{end for}\\
    \STATE $gBest_{i}$=min \{$pBest_{i}$\}
    \STATE \textbf{while} not stop\\
    \STATE ~~ \textbf{for} $i$ =1 to $M$
    \STATE ~~~~~~Updating the velocity and position of particle $i$
    \STATE ~~~~~~Evaluating particle $i$
    \STATE ~~~~~~\textbf{if} fit $({X_{id}}) <$ fit $(pBest_{i})$
    \STATE ~~~~~~~~$pBest_{i} = {X_{id}}$
    \STATE ~~~~~~\textbf{if} fit $(pBest_{i}) <$ fit $(gBest_{i})$
    \STATE ~~~~~~~~$gBest_{i} = pBest_{i}$
    \STATE ~~\textbf{end for}
    \STATE \textbf{end while} \\
    \STATE print $gBest_{i}$
    \STATE \textbf{end}
    \end{algorithmic}
    \end{algorithm}

\subsection{Numerical results for the observable case}

\begin{figure*} [ht]
\centering
\includegraphics[width=7.5cm]{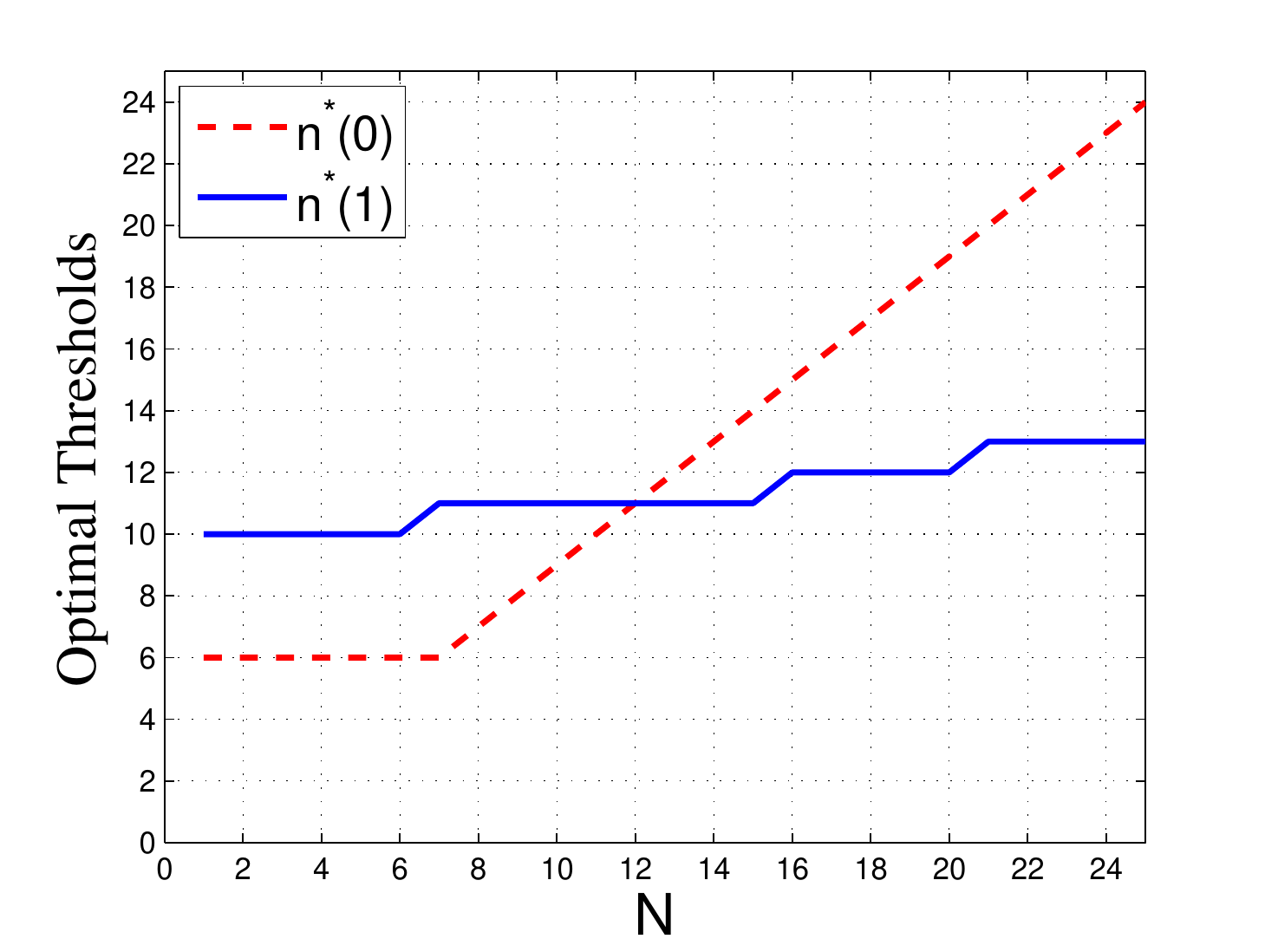}
\caption{Socially optimal thresholds $({n^*}(0), n^*(1))$ with respect to $N$ when $R=15$, $C=1$, $\lambda=5$, $\mu=3$,  $\xi=0.15$, $\theta=5$.}
\label{Fig:6}
\end{figure*}

\begin{figure*} [ht]
\centering
\includegraphics[width=7.5cm]{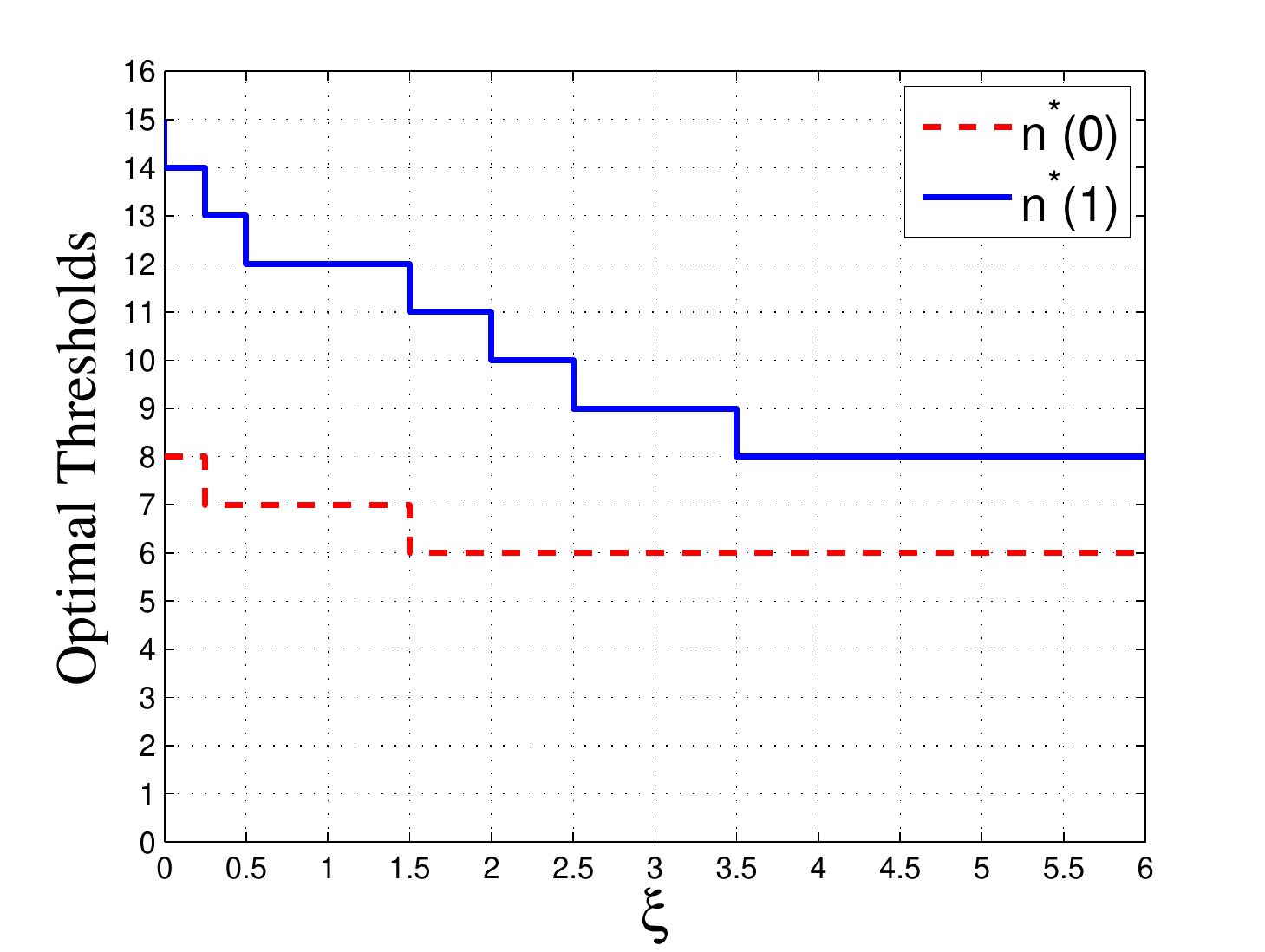}
\caption{Socially optimal thresholds $({n^*}(0), n^*(1))$ with respect to $\xi$ when $R=15$, $C=1$, $\lambda=5$, $\mu=3$,  $\theta=5$, $N=7$.}
\label{Fig:7}
\end{figure*}

Based on a large number of numerical experiments with a series of parameter choices, we conclude that key qualitative properties are independent of the choice of parameters. To illustrate these properties, we present some exemplary results below. First, we explore the trend in changes for the socially optimal thresholds $({n^*}(0),{n^*}(1))$ with respect to $N$ and $\xi$ in Fig \ref{Fig:6} and Fig \ref{Fig:7}, respectively. They illustrate the following phenomena.

1. From Fig. \ref{Fig:6}, we can observe that ${n^*}(0)$ and ${n^*}(1)$ increase with $N$, which illustrates that the social planner wants customers to actively join the system with the growth of $N$.

2. In Fig. \ref{Fig:6}, it is clear ${n^*}(0) \leq {n^*}(1)$ when $N \leq 12$, ${n^*}(0) > {n^*}(1)$ when $N > 12$. The reason for this is that ${U_s}({n^*}(0), n^*(1)) = U_{ob1}^e(n(0),n(1))$ when $N \leq 12$, and ${U_s}({n^*}(0), n^*(1)) = U_{ob3}^e(n(0),n(1))$ when $N > 12$. Specially, ${n^*}(0)=N-1$ when $N>7$.  It indicates that ${n^*}(0)$ is the minimum threshold to ensure server activity. Therefore, when the social planner sets a larger $N$ value, the corresponding ${n^*}(0)$ will be generated.

3. From Fig. \ref{Fig:7}, we can observe that both ${n^*}(0)$ and ${n^*}(1)$ decrease with $\xi$, which illustrates that customers' selfishness does not match the wishes of the social planner. Specially, ${n^*}(0)=N-1$ when $\xi>1.5$, which illustrates that the social planner does not want to accumulate many customers during the vocation.

\begin{figure}[ht]
\centering
\subfigure[]{
\includegraphics[width=7.5cm]{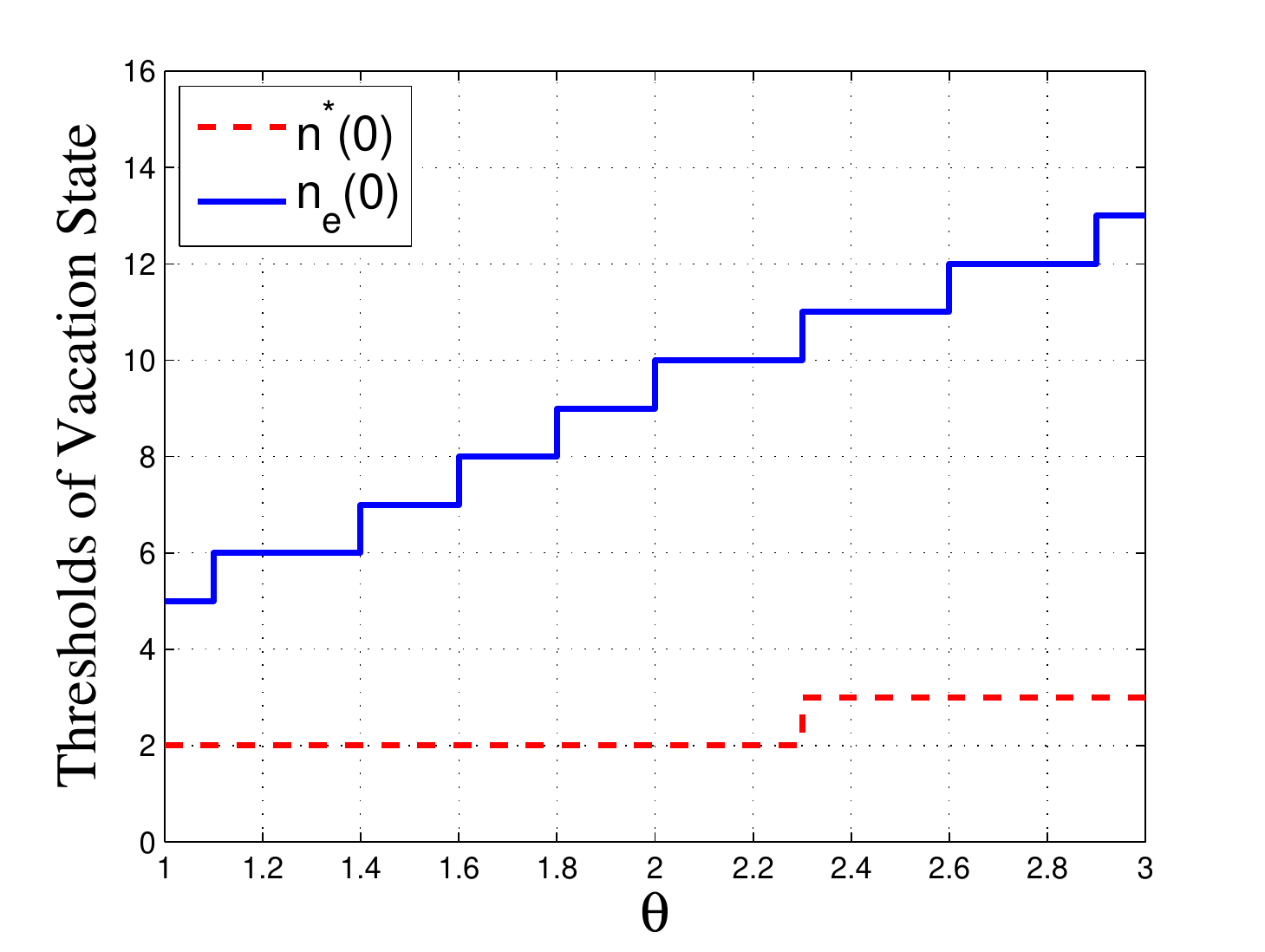}
}
\quad
\subfigure[]{
\includegraphics[width=7.5cm]{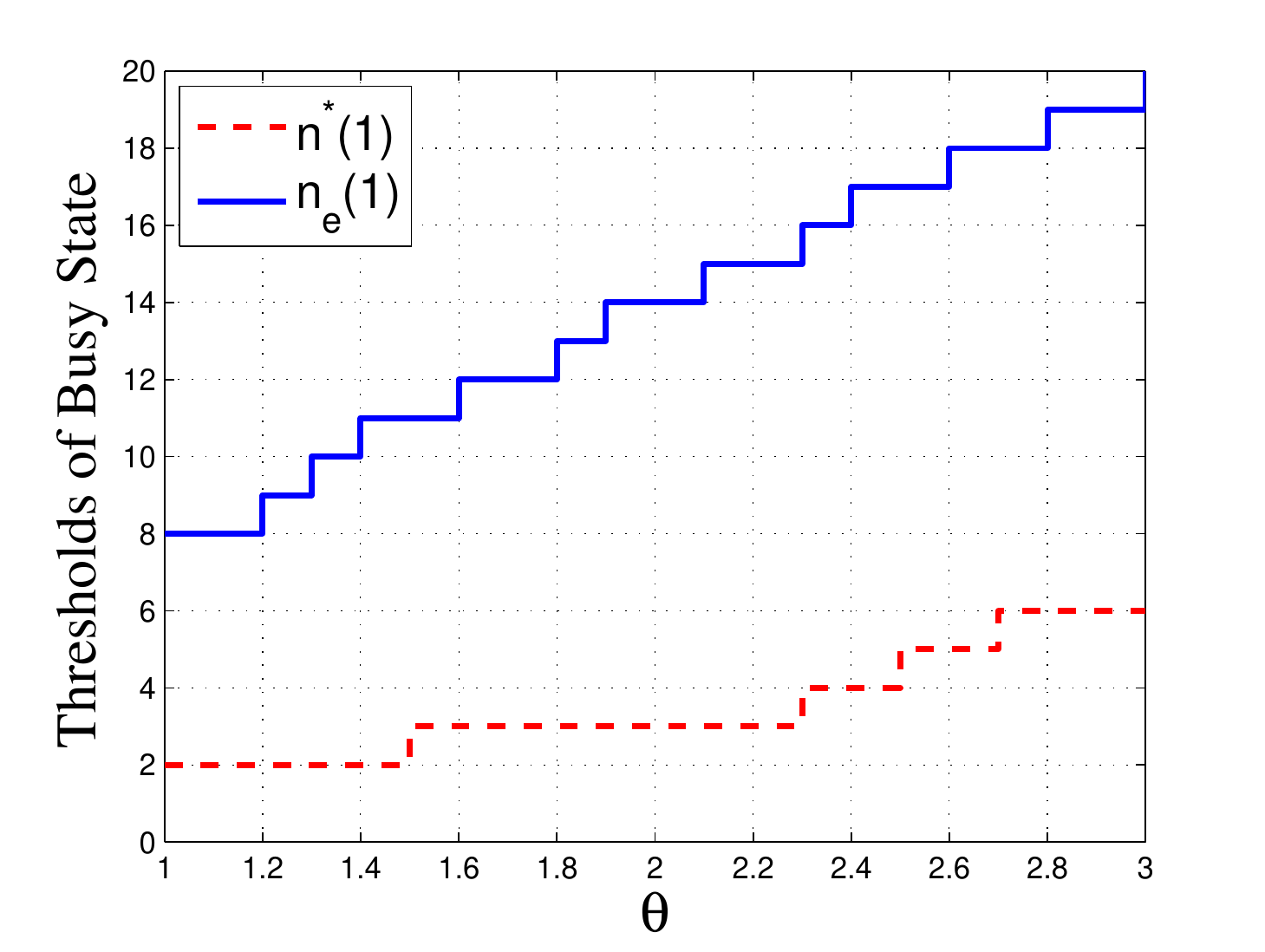}
}
\caption{ { Equilibrium and socially optimal thresholds with respect to $\theta$ when $\lambda=3$, $\mu=5$, $R=15$, $C=1$, $\xi=0.2$, $N=3$.}}
\label{Fig:8}
\end{figure}

\begin{figure}[ht]
\centering
\subfigure[]{
\includegraphics[width=7.5cm]{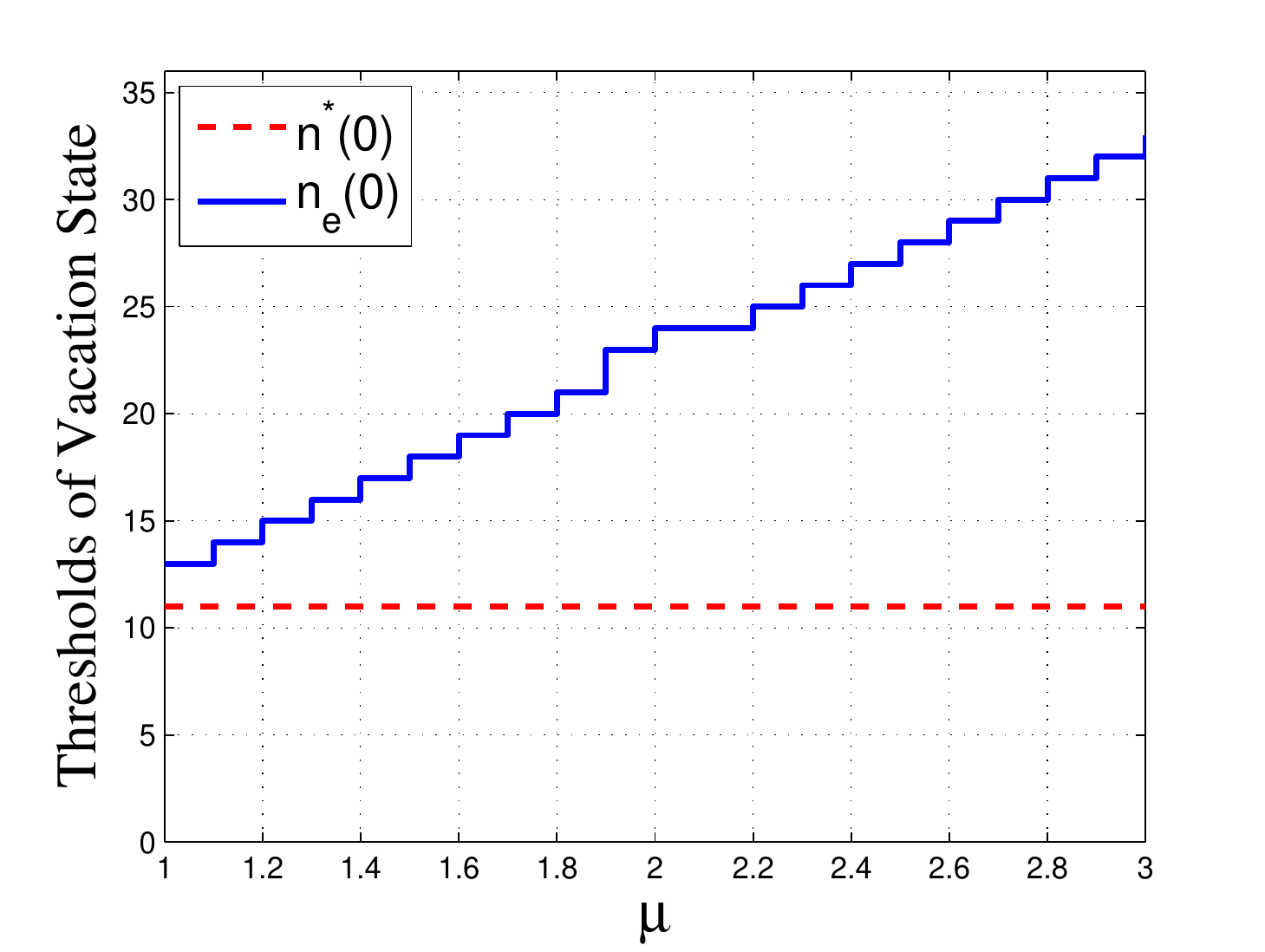}
}
\quad
\subfigure[]{
\includegraphics[width=7.5cm]{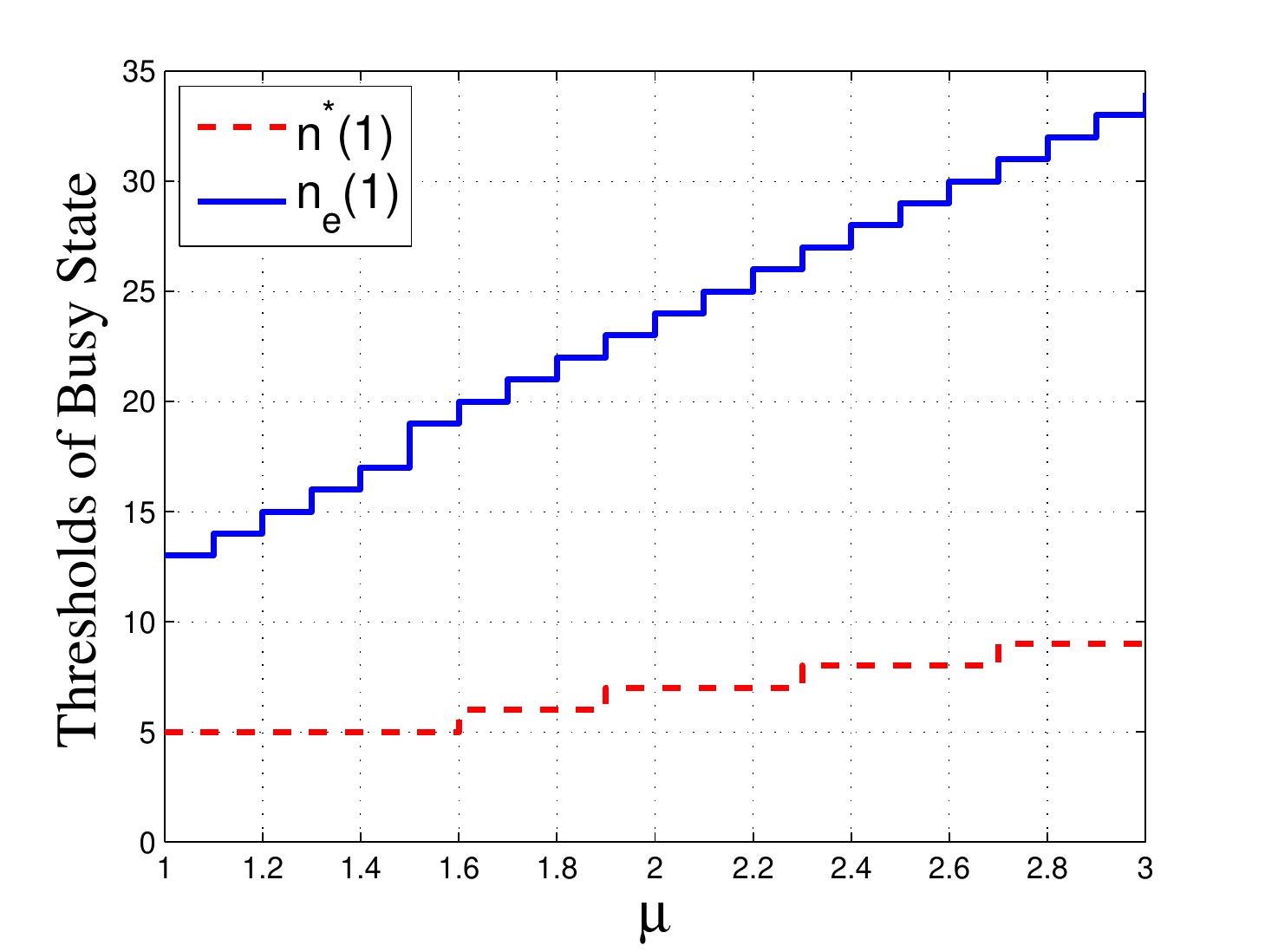}
}
\caption{ Equilibrium and socially optimal thresholds with respect to $\mu$ when $\lambda=3$, $\theta=6$, $R=23$, $C=1$, $\xi=1$, $N=12$.}
\label{Fig:9}
\end{figure}

Next, we explore the trend in changes for ${n^*}(i)$ and $n_{e}(i)$ $(i=0,1)$ with respect to $\theta$ and $\mu$ in Fig \ref{Fig:8} and Fig \ref{Fig:9}, respectively. They reveal the following phenomena.

1. From Fig. \ref{Fig:8} and Fig. \ref{Fig:9}, we can observe that both ${n^*}(i)$ and $n_{e}(i)$ $(i=0,1)$ increase with respect to $\theta$ and $\mu$, respectively. It is obvious that the growth rate of $n_{e}(i)$ $(i=0,1)$ is much faster than the growth rate of ${n^*}(i)$ $(i=0,1)$.

2. $n_{e}(i) > {n^*}(i)$ $(i=0,1)$ always holds as shown in Fig. \ref{Fig:8} and Fig. \ref{Fig:9}, which illustrates that the customers' individual behavior under the stable equilibrium can lead to system congestion more seriously.

\begin{figure}[ht]
\centering
\subfigure[]{
\includegraphics[width=7.5cm]{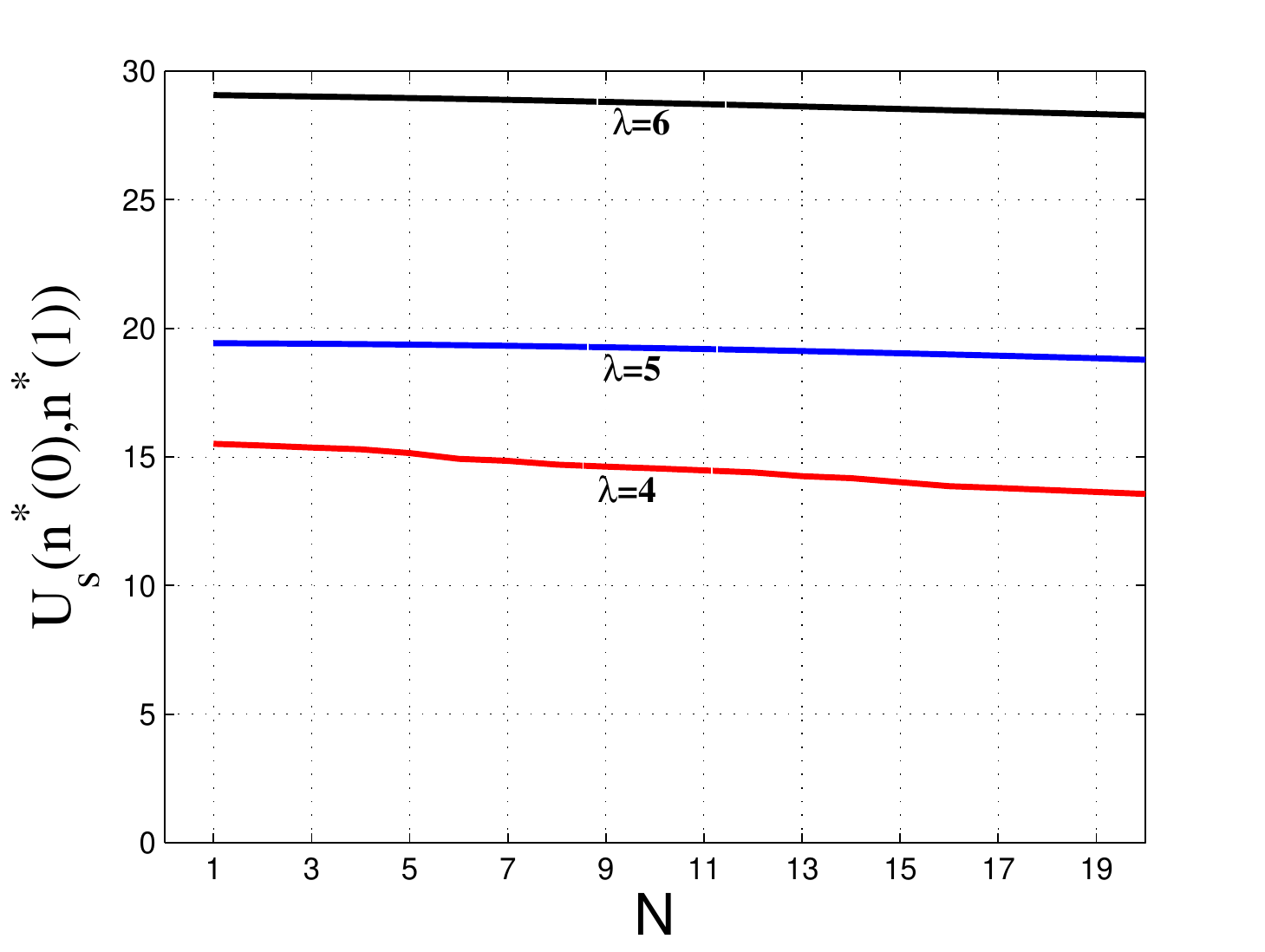}
}
\quad
\subfigure[]{
\includegraphics[width=7.5cm]{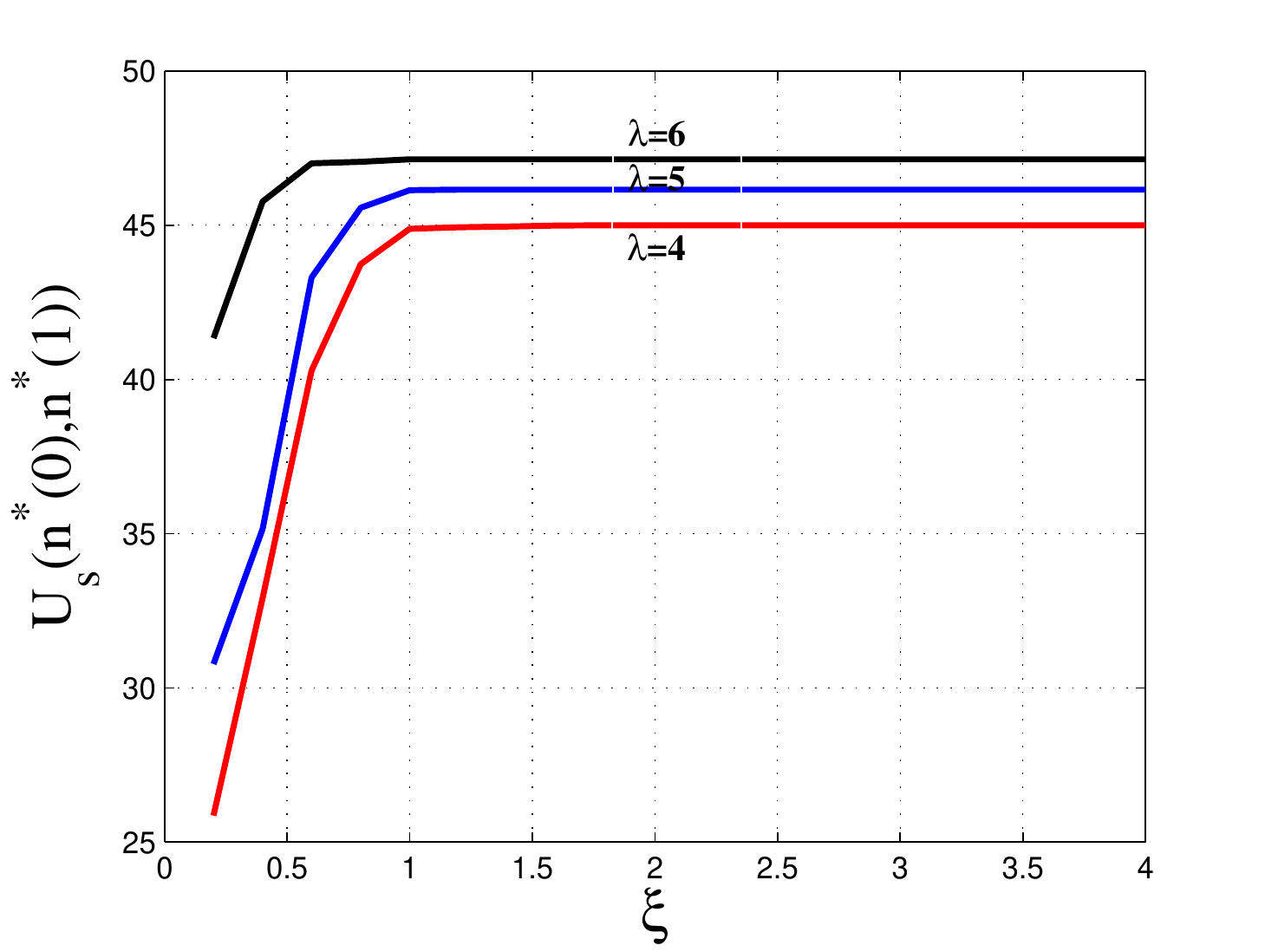}
}
\caption{ (a) Optimal social welfare with $N$ when $\mu=3$, $\xi=0.1$, $\theta=5$, $R=20$, $C=1$; (b) Optimal social welfare with $\xi$ when $\mu=3$, $\theta=5$, $N=6$, $R=20$, $C=1$.}
\label{Fig:10}
\end{figure}

Finally, Fig. \ref{Fig:10} shows that the relationship between the optimal social welfare ${U_s}({n^*}(0),n*(1))$ and $N$,  and the relationship between the optimal social welfare ${U_s}({n^*}(0),n*(1))$ and $\xi$. It reveals the following phenomena.

1. In Fig. \ref{Fig:10} (a), the optimal social welfare ${U_s}({n^*}(0),n*(1))$ decreases with $N$. The reason is that when $N$ becomes larger, more customers will be accumulated, which leads to more waiting costs.

2. In Fig. \ref{Fig:10} (b), the optimal social welfare ${U_s}({n^*}(0),n*(1))$ increases with $\xi$. When $\xi$ becomes bigger, it speeds up the operation of the system, which can then produce the social welfare more effectively. Moreover, when $\xi$ increases to a certain value, the social welfare reaches its maximum and remains stable afterwards.

\subsection{Numerical results for the  unobservable case}

\begin{figure}[ht]
\centering
\subfigure[]{
\includegraphics[width=7.5cm]{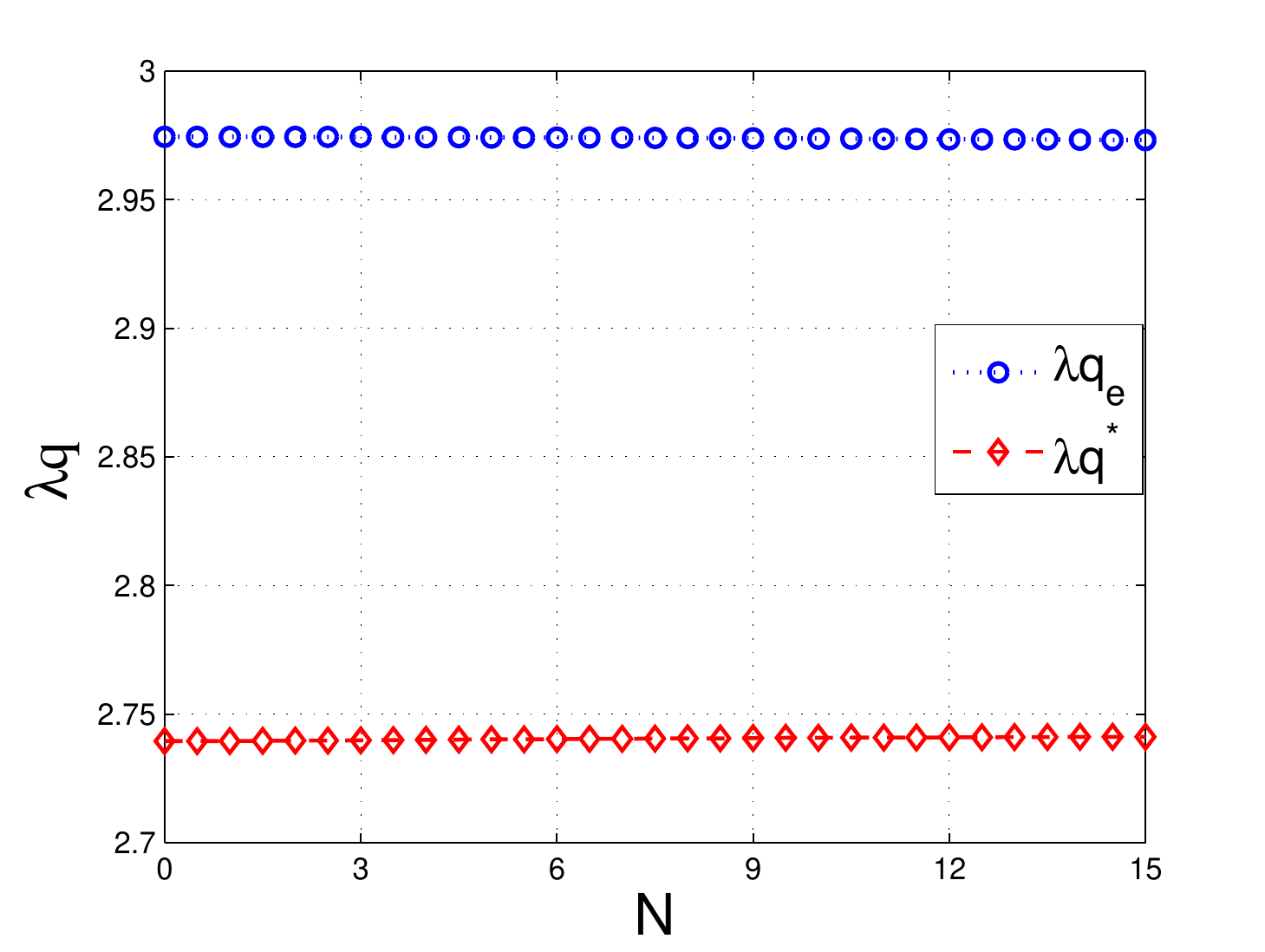}
}
\quad
\subfigure[]{
\includegraphics[width=7.5cm]{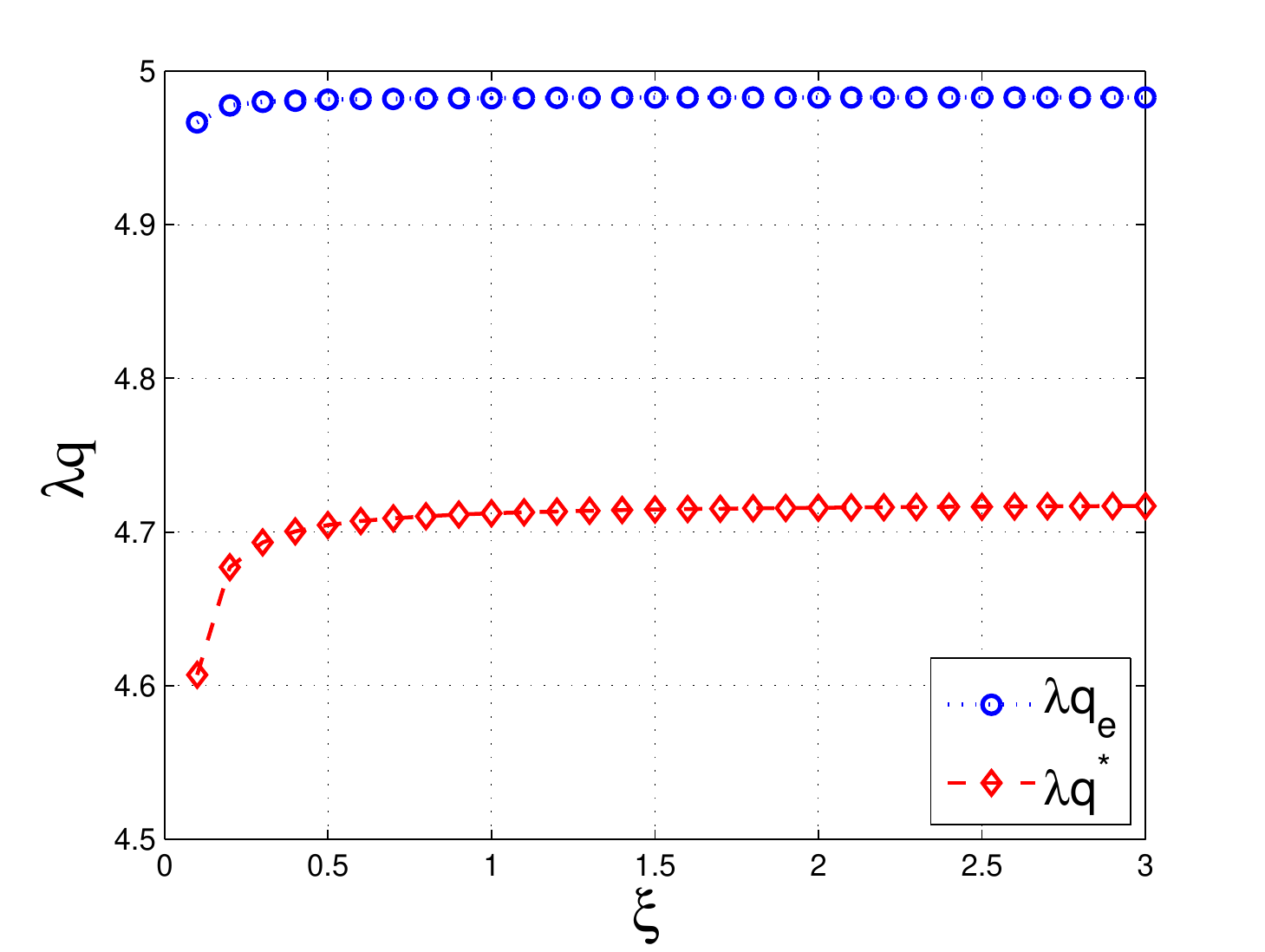}
}
\vfill
\subfigure[]{
\includegraphics[width=7.5cm]{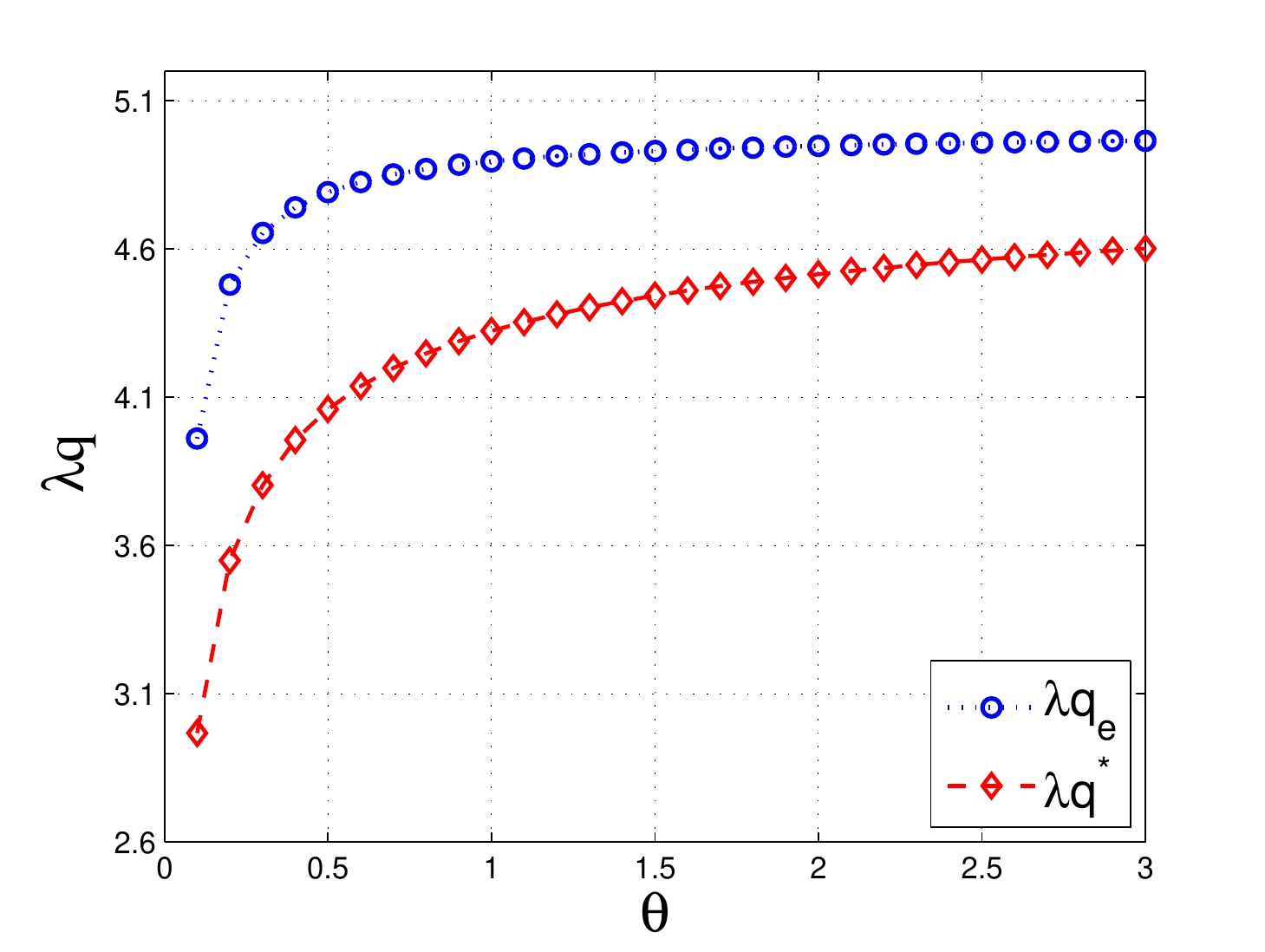}
}
\quad
\subfigure[]{
\includegraphics[width=7.5cm]{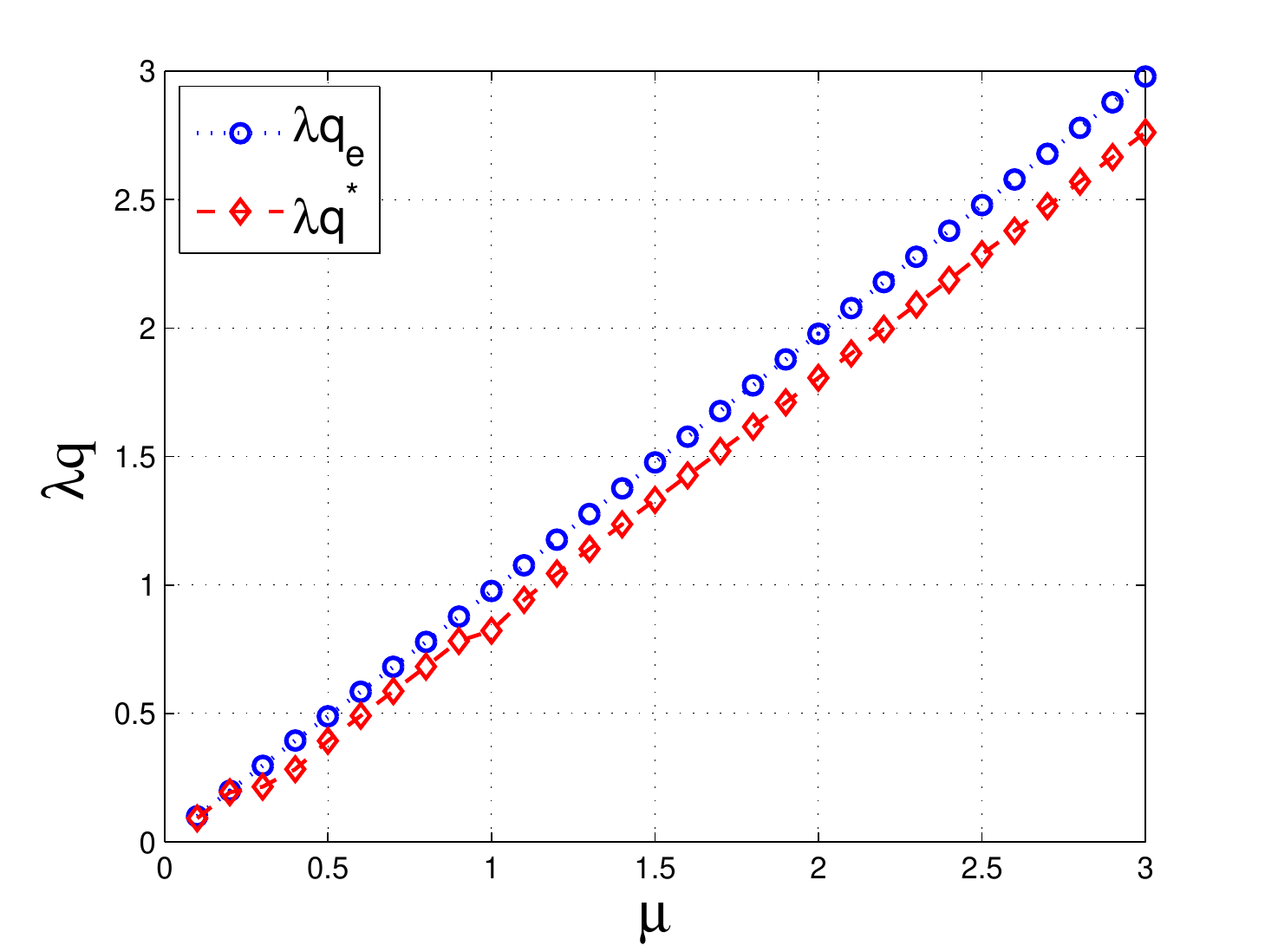}
}
\caption{ Equilibrium and optimal arrival rates with respect to $N$, $\xi$, $\theta$ and $\mu$, respectively. (a) $R=10$, $C=1$, $\lambda=3$, $\mu=3$, $\xi=0.5$, $\theta=5$; (b) $R=20$, $C=1$, $\lambda=5$, $\mu=5$, $\theta=3$, $N=6$ ; (c) $R=10$, $C=1$, $\lambda=5$, $\mu=5$, $\xi=3$, $N=6$; (d) $R=10$, $C=1$, $\lambda=3$, $\xi=3$, $\theta=5$, $N=3$.}
\label{Fig:11}
\end{figure}

In the unobservable case, $\overline {{\lambda _2}} $  is unstable from the equilibrium point of view. Thus, we only study the stable one with $\overline {{\lambda _1}} $, $\overline {{\lambda _3}} $ or $\lambda$ in the following numerical results. First, we explore the impact of the parameters $N$, $\xi$, $\theta$ and $\mu$ on the equilibrium arrival rate $\overline {{\lambda _e}}  = \lambda {q_e}$ and the optimal arrival rate ${\overline \lambda  ^*} = \lambda {q^*}$ in Fig. \ref{Fig:11}, respectively. They illustrate the following phenomena.

1. From Fig. \ref{Fig:11}, we can always see that $\overline {{\lambda _e}}$ (i.e. $\lambda {q_e}$) $\ge {\overline \lambda  ^*}$ (i.e. $\lambda {q^*})$ with $N$, $\xi$, $\theta$ and $\mu$, respectively.

2. In Fig. \ref{Fig:11} (a), $\overline {{\lambda _e}}$ (i.e. $\lambda {q_e}$) decreases with $N$, whereas ${\overline \lambda  ^*}$ (i.e. $\lambda {q^*})$ increases with $N$. The reason for this is that a larger value of $N$ reduces the enthusiasm of customers to join the system. However, a larger value of $N$ makes the social planner encourages more customers to enter the system.

3. From Fig. \ref{Fig:11} (b), (c) and (d), we can observe that the interests of customers and the social planner are coincident with each other with respect to $\xi$, $\theta$ and $\mu$. This is because that shortening the vacation time can reduce waiting costs, and increasing the retrial rate and service rate can accelerate switchover of the system back from vacations.

\begin{figure}[ht]
\centering
\subfigure[]{
\includegraphics[width=7.5cm]{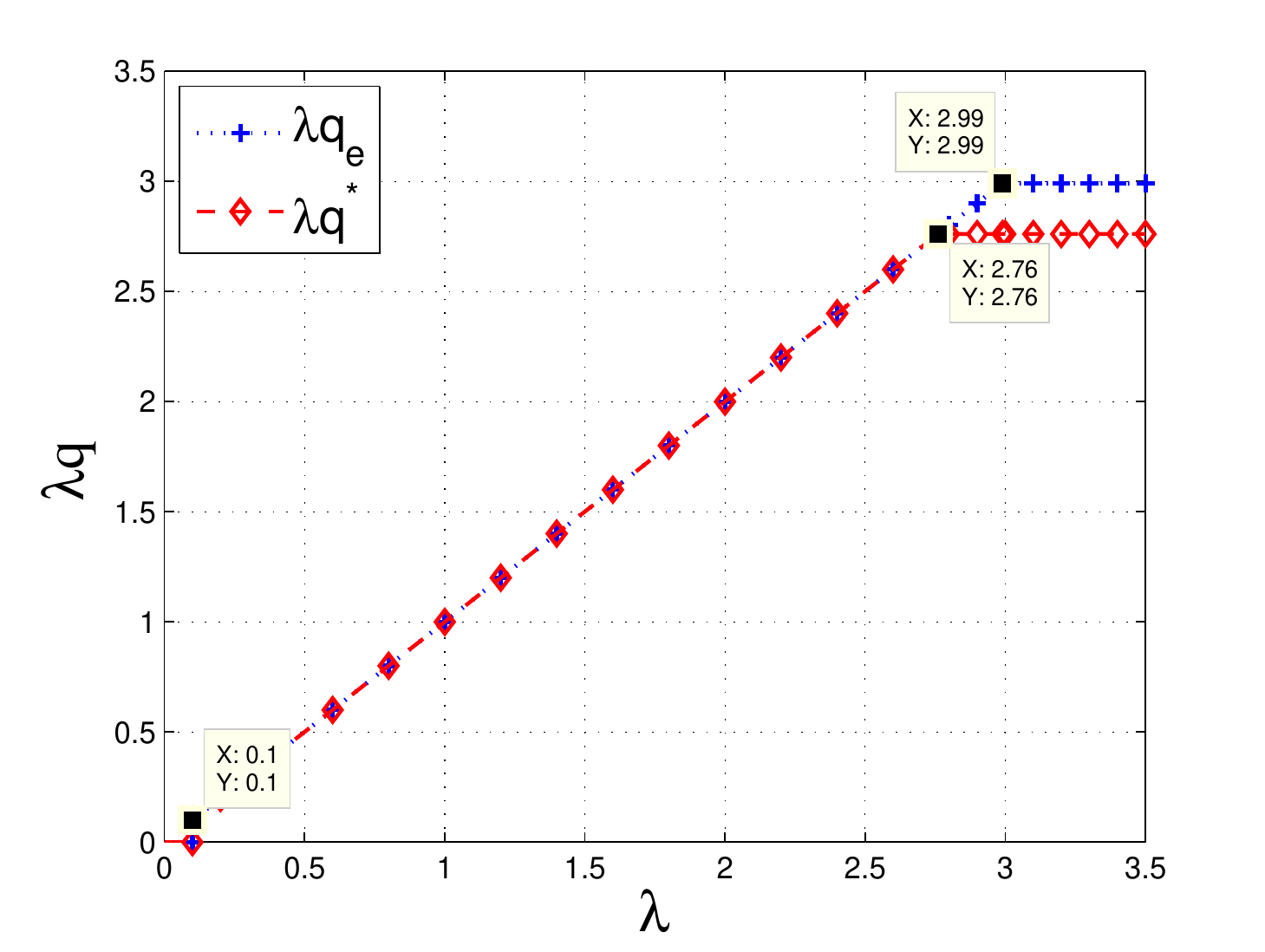}
}
\quad
\subfigure[]{
\includegraphics[width=7.5cm]{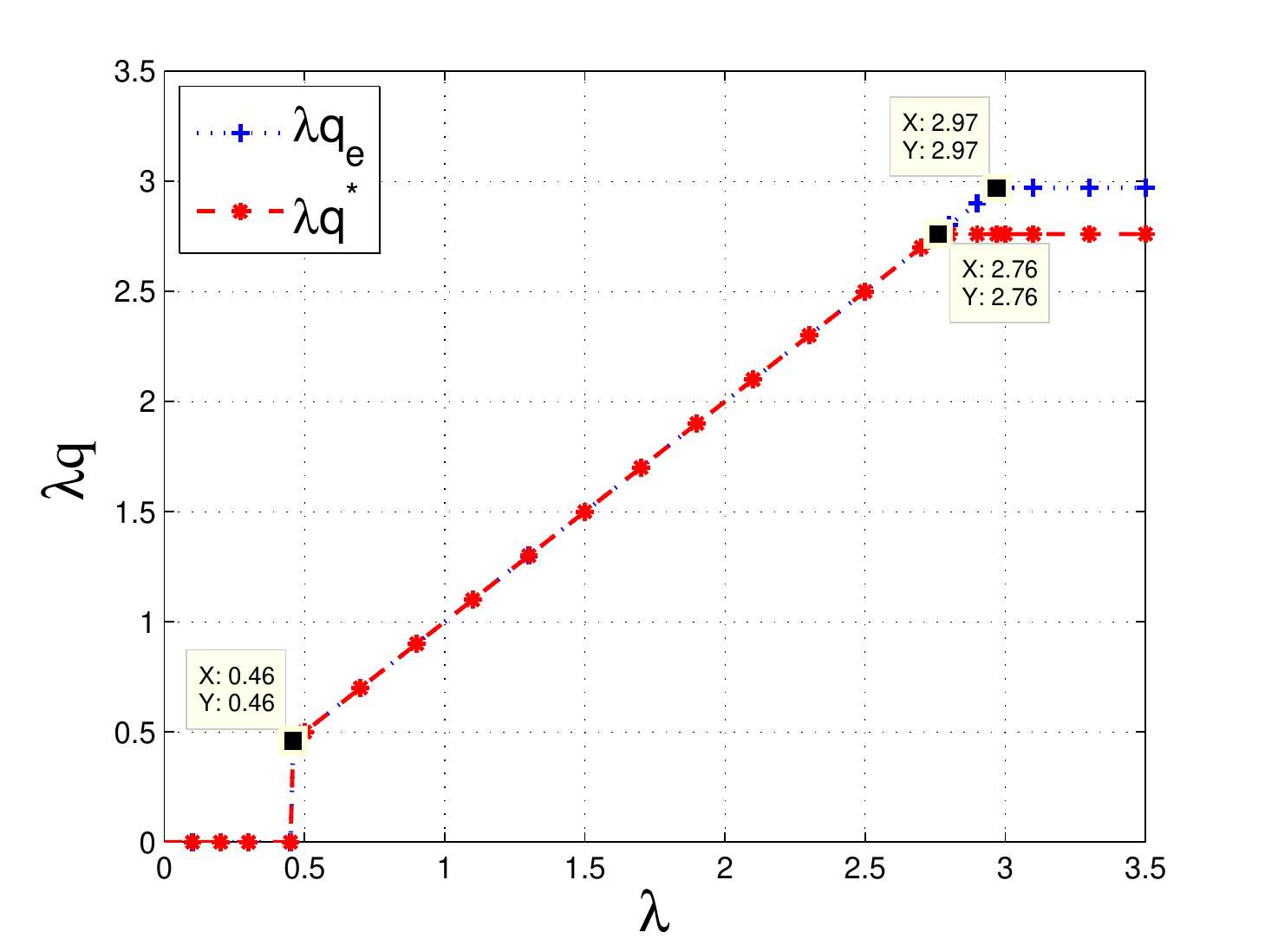}
}
\caption{  Equilibrium and optimal arrival rates for unobservable case with $R=10$, $C=1$, $\mu=3$, $\xi=2$, $\theta=5$; (a) $N=3$; (b) $N=25$.}
\label{Fig:12}
\end{figure}

Moreover, Fig. \ref{Fig:12} shows that $\overline {{\lambda _e}}$ (i.e. $\lambda {q_e}$) $\ge {\overline \lambda  ^*}$ (i.e. $\lambda {q^*})$ as a function of $\lambda$. It is obvious that the value of $N$ has no significance on $\overline {{\lambda _e}}$ and ${\overline \lambda  ^*}$ when $\lambda > \overline {{\lambda _1}} $, but $\overline {{\lambda _1}} $ increases with $N$. Thus, a higher value of $N$ will scare away customers unless $\lambda$ is big enough.

\begin{figure}[ht]
\centering
\subfigure[]{
\includegraphics[width=7.5cm]{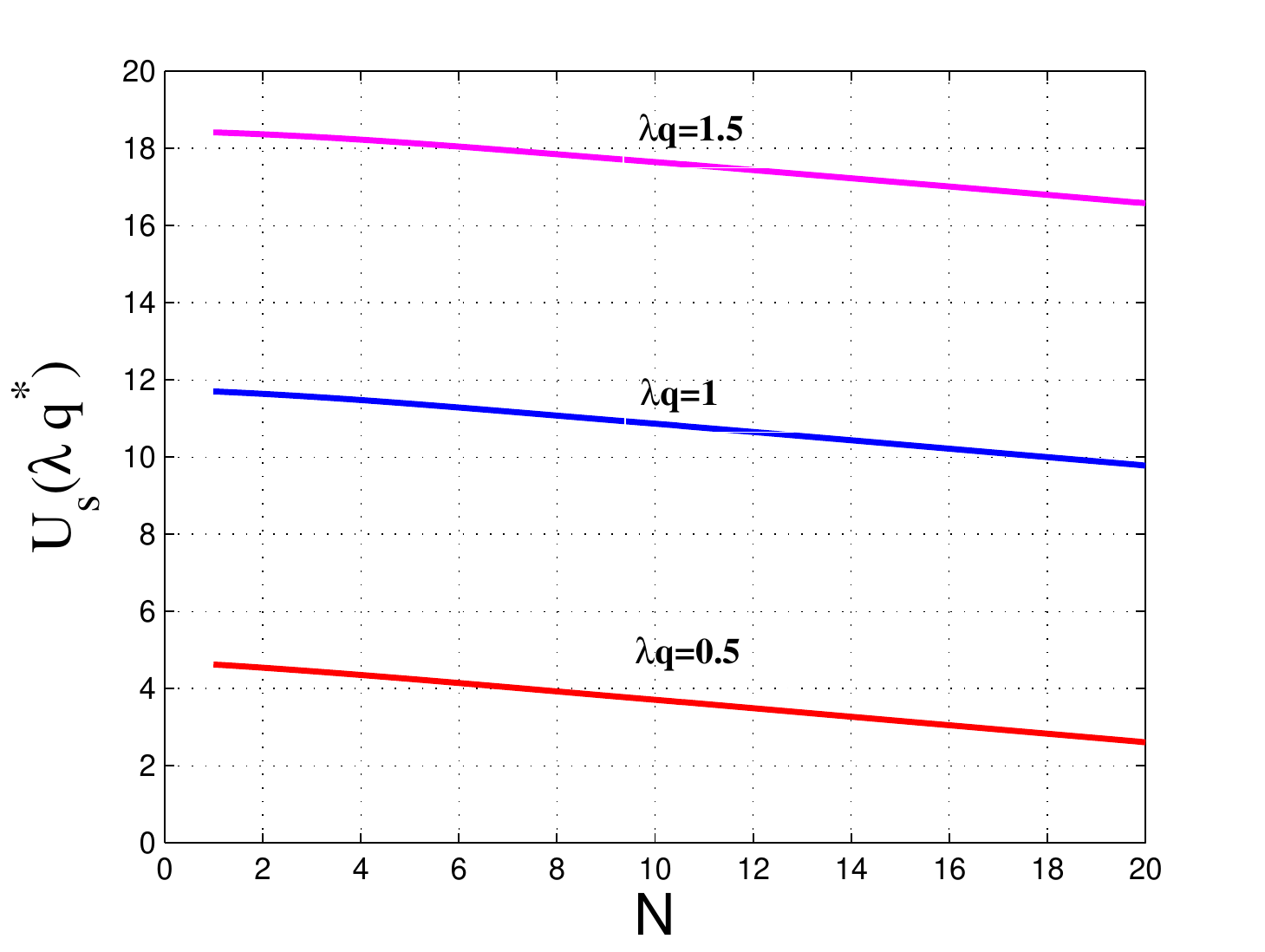}
}
\quad
\subfigure[]{
\includegraphics[width=7.5cm]{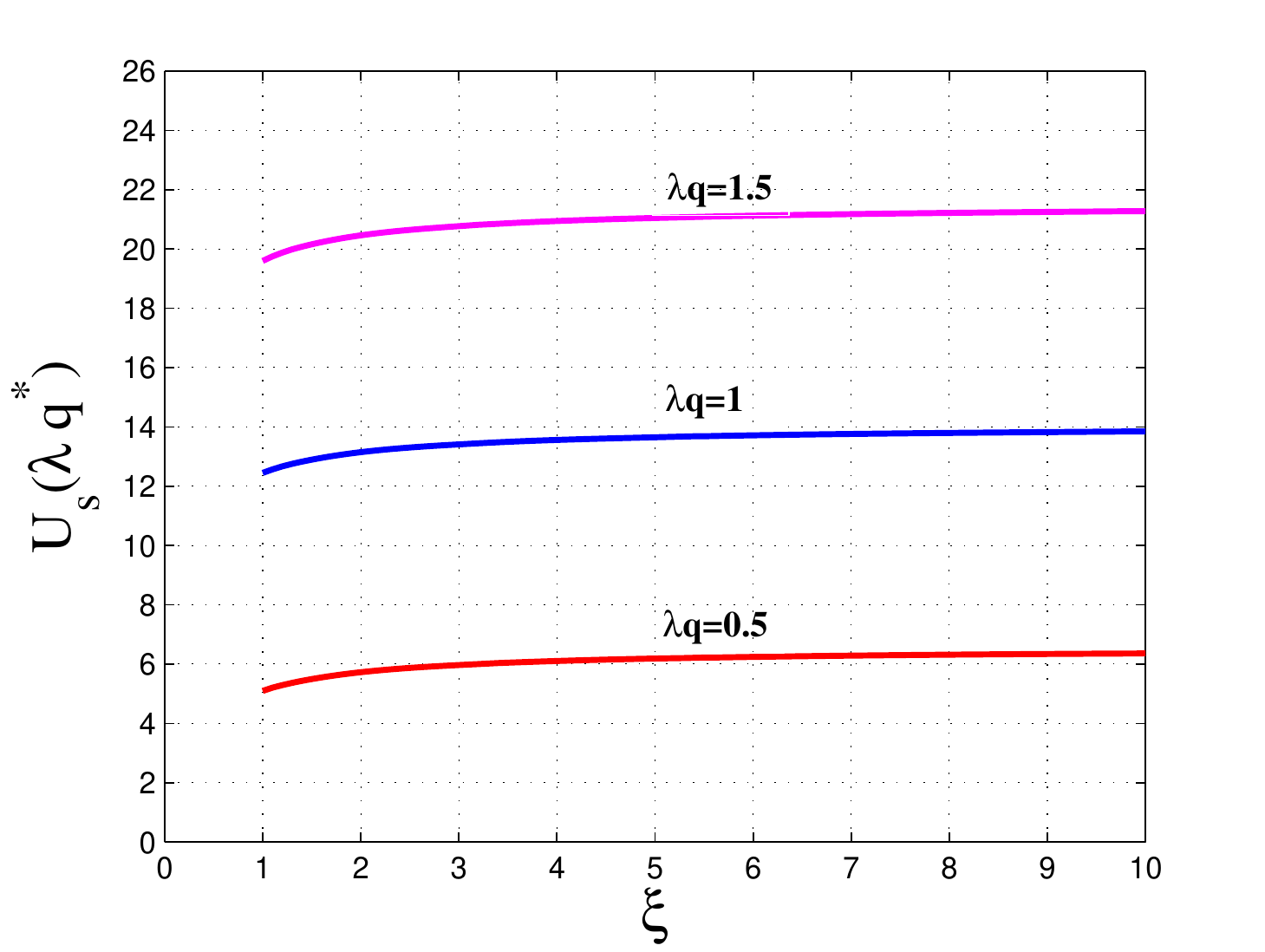}
}
\caption{(a) Optimal social welfare with $N$ when $\mu=3$, $\xi=0.1$, $\theta=5$, $R=20$, $C=1$; (b) Optimal social welfare with $\xi$ when $\mu=3$, $\theta=5$, $N=6$, $R=20$, $C=1$.}
\label{Fig:13}
\end{figure}

If all customers follow the stable equilibrium mixed strategy $\overline {{\lambda _e}} $, then their equilibrium social welfare per time unit ${U_s}(\overline {{\lambda _e}} )$ can be achieved. Similarly, the optimal social welfare per time unit ${U_s}({\overline \lambda  ^*})$ can also be obtained. Fig. \ref{Fig:13} shows the trend in changes for the optimal social welfare per time unit ${U_s}({\overline \lambda  ^*})$ with respect to $N$ and $\xi$, respectively. Obviously, the case shown in Fig. \ref{Fig:13} is similar to that shwon in Fig. \ref{Fig:10}.

\subsection{The role of the information level on the equilibrium social welfare and optimal social welfare}

\begin{figure}[ht]
\centering
\subfigure[]{
\includegraphics[width=7.5cm]{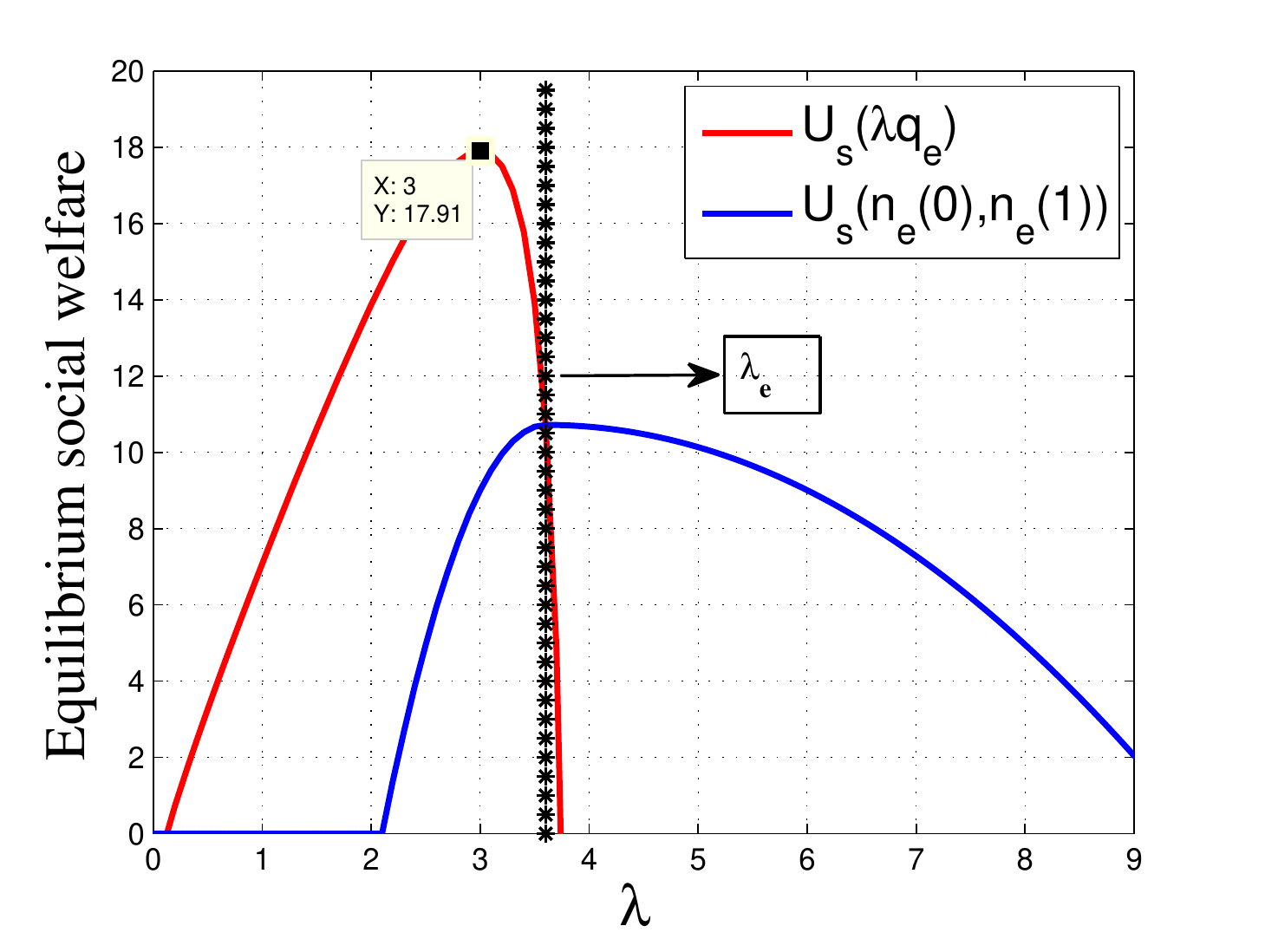}
}
\quad
\subfigure[]{
\includegraphics[width=7.5cm]{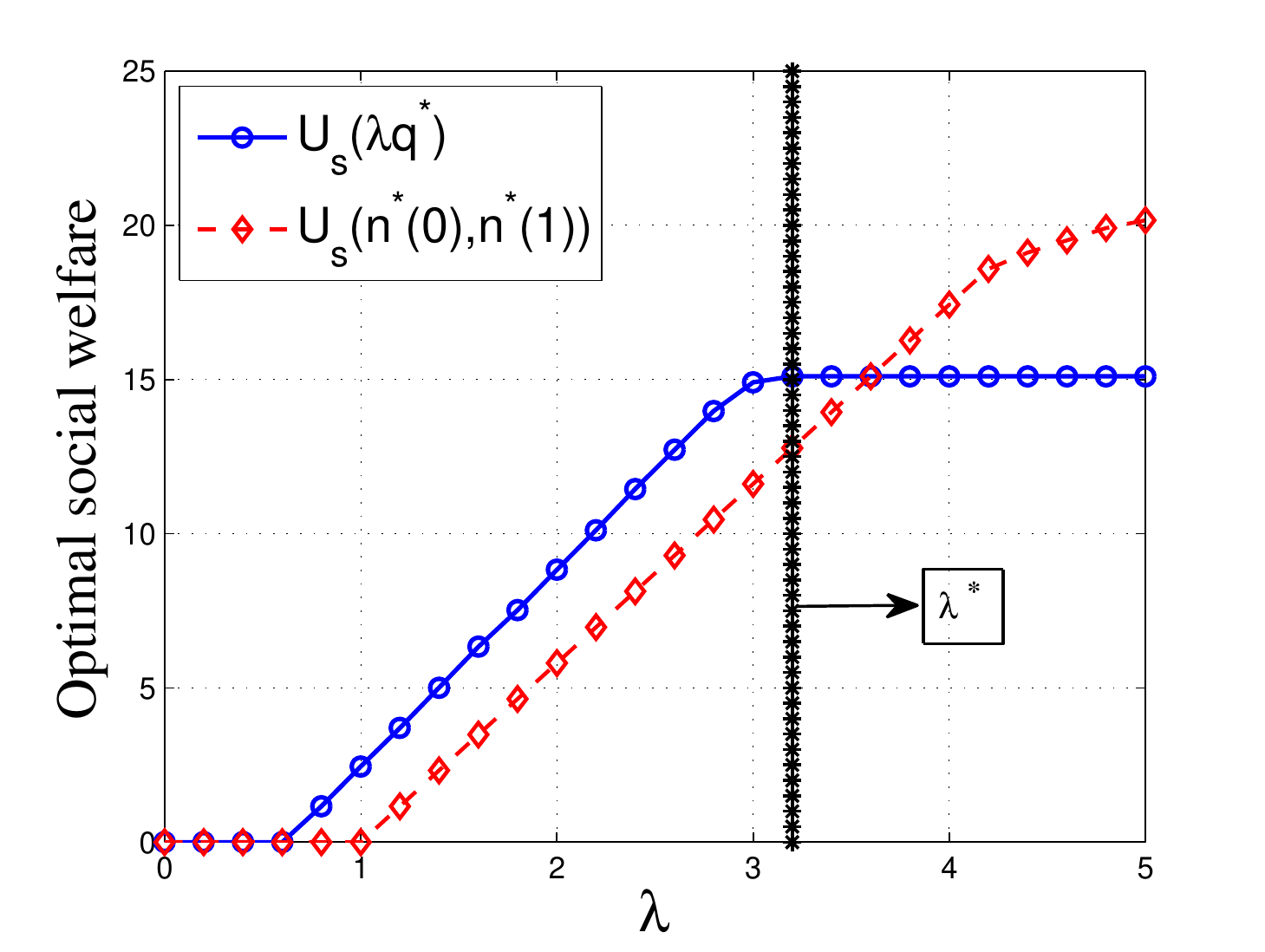}
}
\caption{(a) Comparison of equilibrium social welfare under two information levels when $\mu=3$, $\xi=0.2$, $\theta=5$, $R=20$, $C=1$, $N=6$; (a) Comparison of optimal social welfare under two information levels when $\mu=3$, $\xi=0.2$, $\theta=5$, $R=20$, $C=1$, $N=6$.}
\label{Fig:14}
\end{figure}

An important issue in the strategic customer queuing model is the level of information that the social planner should provide to customers. Fig. \ref{Fig:14} explores the trend in changes for the equilibrium social welfare of customers and optimal social welfare of customers under two levels of information. The properties shown here are presentative since the conclusions made are based on comprehensive numerical experiments with a broad choice of system parameters.

1. In Fig. \ref{Fig:14} (a), it shows that the equilibrium strategy of customers is balking when $\lambda$ is small, since smaller $\lambda$ does not intend to activate the system. In this case, both ${U_s}({n_e}(0),{n_e}(1))$ and ${U_s}(\overline {{\lambda _e}} )$ (i.e. ${U_s}(\lambda {q_e})$) have a zero-increase-decrease trend with $\lambda$. However, ${U_s}(\overline {{\lambda _e}} ) \ge {U_s}({n_e}(0),{n_e}(1))$ when $\lambda  < {\lambda _e}$, since the arrival rate is small and the number of people in the system is small. In this case, hiding the system information from customers helps to increase the number of customers entering the system, thereby increasing social welfare. Similarly, ${U_s}(\overline {{\lambda _e}} ) < {U_s}({n_e}(0),{n_e}(1))$ when $\lambda  > {\lambda _e}$, which shows that disclosing the system information can help reduce the system congestion, thus reducing waiting costs.

2. In Fig. \ref{Fig:14} (b),  it shows that the socially optimal mixed strategy of customers is balking when $\lambda$ is small. Both ${U_s}({n^*}(0),{n^*}(1))$ and ${U_s}({\overline \lambda  ^*})$ (i.e. ${U_s}(\lambda {q^*})$) keep growth, ${U_s}({\overline \lambda  ^*})$ eventually becomes a constant. Similar to Fig. \ref{Fig:14} (a), Fig. \ref{Fig:14} (b) also shows that ${U_s}({\overline \lambda  ^*}) \ge {U_s}({n^*}(0),{n^*}(1))$ when $\lambda  < {\lambda ^*}$, ${U_s}({\overline \lambda  ^*}) < {U_s}({n^*}(0),{n^*}(1))$ when $\lambda  > {\lambda ^*}$. This also shows that the information level of the system has a serious impact on the social welfare. Therefore, the social planner should choose the strategy consistent with the system designer to achieve social optimum.

\section{Conclusions and further research}
\label{sec:6}
In this paper, we studied equilibrium strategies and optimal balking strategies of customers in a constant retrial queue with multiple vacations and the $N$-policy under two information levels (observable case and unobservable case), respectively. For each type of information levels, we determined equilibrium strategies and optimal balking strategies of customers, and the social welfare. For the observable case, in order to ensure that the server can be reactivated, we obtained that the optimal balking threshold of customers in the vacation state must be greater than the optimal threshold in busy state or greater than $N$. Therefore, there are three different queuing cases for the observable case, and we obtained the corresponding stationary distributions for the three queuing cases, and determined the equilibrium social welfare per time unit. For the unobservable case, we obtained the positive equilibrium arrival rate and optimal arrival rate, which are unique. In Section \ref{sec:5}, we explored the previous theoretical results through numerical experiments. However, due to the complexity of the involved equations, explicit expressions for the equilibrium balking thresholds of customers, socially optimal balking thresholds and optimal social welfare are not available in general. Hence, we use Particle Swarm Optimization (PSO) algorithm to solve the complex analytic characteristics. The numerical optimal solution $({n^*}(1),{n^*}(2))$ and optimal social welfare ${U_s}({n^*}(0),{n^*}(1))$ and ${U_s}({\overline \lambda  ^*})$ are obtained by PSO algorithm. By comparing the numerical results of the two information levels, we obtained that the customers' behavior under the stable equilibrium makes the system more congested than that under the socially optimal one, and whether the system information should be disclosed to customers depends on how to maintain the growth of the social welfare (i.e., potential demand arrivals). Obviously, in order to maximize the social welfare, which factor determines the level of information disclosure and when to disclose the system information to customers are also crucial for the server or social planner. Fortunately, this paper achieved this goal. In the future, it is necessary for us to consider the almost observable case of this model, i.e., the state of the server can be observed, but the number of customers in the orbit cannot be observed.

\section*{Acknowledgements}
This work was supported in partial by The National Natural Science Foundation of China (No. 61773014), the program of China Scholarships Council (No. 201906840070), and by The Natural Sciences and Engineering Research Council of Canada (NSERC).

\section*{Disclosure statement}
None of the authors have any competing interests in the manuscript.


\end{document}